%% file: Main.tex
\documentclass{elsarticle}

\input{tex/sty}

\definecolor{rev}{gray}{0.4} 
\title{BATTPOWER Application: Large-Scale Integration of EVs in an Active Distribution Grid ---A Norwegian Case Study}

\author[a]{Salman Zaferanlouei}
\author[b]{Venkatachalam Lakshmanan}
\author[a]{Sigurd Bjarghov}
\author[a]{Hossein Farahmand}
\author[a]{Magnus Korpås}

\address[a]{Department of Electric Power Engineering NTNU Trondheim, Norway}
\address[b]{SINTEF Energy Research, Sem Sælands vei 11, 7034 Trondheim, Norway.}

\begin{document}
	\newgeometry{left=2.2cm,right=2.5cm,top=3cm,bottom=3cm}
\input{tex/abstract}


\maketitle
\input{tex/intro.tex}
\input{tex/problemFormulate.tex}
\input{tex/inpData.tex}

\input{tex/methodology.tex}
\input{tex/result.tex}
\input{tex/conclusion.tex}
\input{tex/acknow.tex}
{\footnotesize
\bibliographystyle{ieeetr}
\bibliography{refs} } 
\appendix
 \input{tex/appA.tex}
 \input{tex/appB}
 \input{tex/appC}
 \input{tex/appD}

 \input{tex/appE}

 \input{tex/appF}

\end{document}

%% file: tex/sty.tex
\usepackage{color}
\usepackage{lineno}
\usepackage{pdflscape}
\usepackage{amsmath}
\usepackage{stackengine}
\usepackage{todonotes}
\usepackage{pdflscape}
\usepackage{rotating}
\usepackage{booktabs}
\usepackage[hyphens]{url}
\usepackage{hyperref}
\usepackage{color}
\definecolor{darkblue}{rgb}{0,0.2,0.4}

\definecolor{rev}{gray}{0.4} 
\usepackage{natbib} 
\biboptions{sort&compress}

\usepackage{colortbl}
\newcolumntype{L}{>{\centering\arraybackslash}m{3cm}}

\usepackage{multirow}

\usepackage{caption}
\usepackage{float}
\usepackage[pass]{geometry}

\usepackage[perpage]{footmisc}

\usepackage[linesnumbered,ruled,vlined]{algorithm2e}

\SetCommentSty{mycommfont}

\usepackage{bookmark}
\usepackage{lipsum}
\usepackage{tabulary,booktabs}

{}
{}
{}
\usepackage{amsmath,amsfonts,amssymb}
\usepackage{mathtools,amsmath}

\usepackage{threeparttable}
\usepackage{nomencl}
\usepackage{makeidx}

\makenomenclature
\usepackage{url}
\usepackage{etoolbox}
 \renewcommand\nomgroup[1]{%
  \item[\bfseries
  \ifstrequal{#1}{A}{Abbreviations}{%
  \ifstrequal{#1}{B}{General}{%
  \ifstrequal{#1}{C}{{BATTPOWER Input Data, (more detailed description in \ref{appendixC})}}{%
  \ifstrequal{#1}{D}{BATTPOWER Parameters}{%
  \ifstrequal{#1}{E}{BATTPOWER Variables}{%
  \ifstrequal{#1}{F}{BATTPOWER Indices}{%
  \ifstrequal{#1}{G}{Algorithm \ref{Alg1}}{}}}}}}}%
 ]}

\DeclarePairedDelimiter\abs{\lvert}{\rvert}%

\usepackage{tabularx,booktabs}

\newcolumntype{C}{>{\centering\arraybackslash}X} 
 \newcolumntype{b}{>{\centering\arraybackslash\hsize=2.3\hsize}X}
\newcolumntype{s}{>{\centering\arraybackslash\hsize=.45\hsize}X}
\newcolumntype{m}{>{\centering\arraybackslash\hsize=.9\hsize}X}

\usepackage{eurosym}

\usepackage{xparse}
\usepackage{hyperref}
\usepackage{cleveref}

\newcounter{starredsection}%

\makeatletter
\let\latex@@section\section

\AtBeginDocument{%
  \RenewDocumentCommand{\section}{som}{%
    \IfBooleanTF{#1}{%
      \refstepcounter{starredsection}%
      \latex@@section*{#3}%
    }{%
      \IfValueTF{#2}{%
        \latex@@section[#2]{#3}%
      }{%
        \latex@@section{#3}%
      }%
    }%
  }%
}
\makeatother

\usepackage{rotating}
\usepackage{graphicx}
\usepackage[pass]{geometry}

\usepackage{appendix}

\usepackage{afterpage}
\usepackage{amsmath}

%% file: tex/abstract.tex
\begin{abstract}
With the considerable increase of Distributed Energy Resources (DER), the reliable and cost-effective operation of distribution grids becomes challenging. The efficient operation relies on computationally dependable and tractable optimisation solvers, which may handle: 1) non-linear AC power flow constraints, and 2) time-linking variables and constraints and objectives of DER, over the operational horizon. In this paper, we introduce an application of a high-performance MultiPeriod AC Optimal Power Flow (MPOPF) solver, called ``BATTPOWER", to simulate active distribution grids for a near-future scenario. A large-scale Norwegian distribution grid along with a large population of Electric Vehicles (EV) are here taken as the case-study. We suggest and analyse three operational strategies (in terms of control of charge scheduling fleet of EV) for the Distribution System Operator (DSO): (a) uncoordinated/dumb charge scheduling, (b) coordinated charge scheduling with the objective of energy cost-minimisation without operational constraints of the grid, and (c) coordinated charge scheduling with the objective of energy cost-minimisation along with the operational constraints of the grid. The results demonstrate that the uncoordinated charging would lead to: 1) overloading of lines and transformers when the share of EVs is above 20\%, and 2) higher operational costs than the proposed control strategies of (b) and (c). In strategy (b) operational line/transformer limits are violated when the populations of EVs are growing above 36\%. This implies that current market design must be altered to allow active control of a large proportion of DERs within grid operational limits to achieve cost minimization at system level. To our knowledge, the work presented in this paper is the first ever attempt to do a comprehensive analysis of the impact of EV charging demand on a real Norwegian distribution grid. Moreover, the inference of the analysis says that the Norwegian distribution networks are more prone to congestion problems than the voltage problems for the EV demand which includes a smart charging scheme accounting for grid conditions.
\end{abstract}

%% file: tex/intro.tex
\section{Introduction} \label{intro}
In recent years, a large increase in electric vehicle (EV) sales has been observed due to decreasing battery prices, larger production volumes and climate policies \cite{EVoutlook}. With the highest EV share in the world, Norway stands out in particular with more than 15\% of the total car park being full electric, numbering more than 400,000 cars \cite{elbilstatistikk}. With a market share of 55\% of all new sales, the growth is expected to continue. The growing additional demand of EV charging is predicted to create congestion in distribution grids, and the Norwegian regulator estimates that 1.2 billion \euro \ can be saved by smart coordination of EV charging \cite{nveelbil}.\\
Currently, most EV chargers start charging at the nominal charging power until the EV battery is full or until a set point has been reached. In the future, we assume that EV charging can be controlled according to the wishes of the EV owner (which to some extent already exists). Many potential charging schemes are suggested in the literature, such as charging ``queues", ``bandwidth sharing" and price signals. However, these approaches often ignore power system aspects or strongly simplify them. By considering EV charging as part of a MultiPeriod AC Optimal Power Flow (MPOPF) problem, optimal charging schedules considering system feasibility and voltages can be achieved, while still minimising costs. 

\subsection{MultiPeriod AC Optimal Power Flow (MPOPF)} \label{MPOPF}
The optimal power flow (OPF) is a non-linear, non-convex problem, introduced in the 60s \cite{carpentier1962contribution}. Depending on the technical applications and operational dimensions, the OPF problem may be evolved to various versions such as the MultiPeriod AC Optimal Power Flow (MPOPF) \cite{chandy_simple_2010}, and may become intractable and computationally hard due to nonlinearities and poor scalability.\\

Since the non-linear ACOPF problems require non-linear solvers to be called, several Non-Linear Programming (NLP) solvers primarily developed based on Interior Point (IP) methods are used to solve MPOPF problems, such as MIPS \cite{wang_computational_2007}, IPOPT \cite{wachter_implementation_2006}, KNITRO \cite{byrd_knitro_2006}, and recently BELTISTOS \cite{kourounis_towards_2018}. Of the mentioned, only BELTISTOS is tailored to solve MPOPF problems. An extensive review of both MPOPF problem formulations and solution methods can be found in \cite{sperstad_optimal_2016, sperstad_energy_2019}. The authors recommend the highly relevant work in Ref. \cite{moghadasi_optimal_2016}, in which a solution to MPOPF problems has been suggested and tested for online implementation purposes. However, the largest examined case study, for the purpose of scalability test, is the 119-bus distribution test system with four storage devices.\\
It is well-known from the literature \cite{capitanescu_experiments_2013, castillo2013computational} that the solution of linear Karush–Kuhn–Tucker (KKT) systems and calculation of gradients are the two most computationally expensive aspects in solving a MPOPF problem. Thus, we proposed a fast solver to exploit the sparsity of a MPOPF structure (both KKT systems and gradients) and to speed up the solution \cite{zaferanlouei_battpower_2021}. 

\subsection{Centralised Optimal Charge Scheduling of EVs} \label{EVCharge}
A substantial amount of research has been conducted to solve the optimal scheduling of EV charging in distribution grids. Reference \cite{sojoudi_optimal_2011} formulated an optimisation framework by assigning price-elastic load to EV charging and considering optimal power flow as the balance constraint and cost of production as the objective function. The optimisation model is applied to a small-scale IEEE 14-bus system. A Smart Load Management (SLM) system is proposed by \cite{masoum_smart_2011} for the coordination of large-scale EV chargers in distribution feeders. The SLM approach is tested for a 1200 bus test system consisting of low-voltage residential networks. Two EV charging controller methodologies, a local EV charger controller and a centralised EV charger controller, are introduced and compared by  \cite{richardson_local_2012}. They suggested that although the network and communication infrastructure needed to implement the local control method would be far less than that of the centralised control case, the centralised controller gives a more reliable operational outcome in the case of high EV penetration. Reference \cite{chen_optimal_2012} suggested a joint optimal power flow and EV charging framework that considers an OPF problem with EV charging over time. This nested optimisation problem is solved through a decomposition approach which has lower computational complexity than that of centralised interior point solvers. The approach is implemented on a IEEE 14-bus. Reference \cite{oconnell_rolling_2014} proposed an unbalanced three-phase multiperiod AC optimal power flow optimisation problem which allocates individual variables to each EV and controls the charging rate and times of charge of EVs over a 24-hour time horizon. The cost function is to minimise the total cost subject to operational constraints. The proposed formulation is solved through an NLP solver of MATLAB, called FMINCON. Moreover, it is applied on a 85 bus test case, 74 single phase, 11 three-phase case studies. However, there is no discussion as to whether the proposed algorithm is fast enough to be scalable. Reference \cite{franco_mixed-integer_2015} proposed a mixed-integer linear programming (MILP) model for EV charging coordination in an unbalance distribution network. The proposed method takes into account the distributed generators and operational constraints. The linear proposed model is solved using commercial MILP solvers. The proposed model is tested on a 394-bus distribution system. However, there is no discussion over the computational complexity of the proposed problem in the paper. Reference \cite{benetti_real-time_2015} proposed a method for the real-time management of EV charging procedures such that it flattens the peak load, increases the number of rechargeable EVs, and activates the network operational constraints. The approach integrates: 1) the scheduling algorithm, 2) power flow equations, and 3) operational constraints. Simulations are conducted on a real medium-size Italian electricity distribution grid. The size of the simulated grid is unknown in the paper. There is no discussion of the computational complexity of the proposed method in the paper. Reference \cite{de_hoog_optimal_2015} studied a fast receding horizon optimisation problem by linearising voltage drop in the network. Two objectives have been considered: 1) maximisation of total EVs charging in the network, and 2) cost of charging. Higher efficiency to exploit the existing lines/transformers in the distribution network is observed in the simulation results of the proposed method. Reference \cite{shao_layered_2015} proposed a layered distributed charging load for controlled charging of EVs based on Lagrangian relaxation and auxiliary problem principle. The proposed method is suitable for large populations of EVs and gains an advantage in reducing generation cost. However, load flow equations are not considered as a part of the formulation. Reference \cite{wang_integrated_2016} proposed a two-stage energy exchange planning strategy for a multi-micro-grid system incorporating EVs as storage devices. The proposed method brings down the electricity cost and prevents frequent transition between charge and discharge modes. Reference \cite{quiros-tortos_control_2016} proposed: 1) a centralised control algorithm, which uses limited data to manage EV charging stations to mitigate grid operational constraints, and 2) OPF-based method. The first method is implemented on two real and large distribution grid cases with 351 and 428 customers. The proposed control algorithm works based on the selection of time of charge and considers the transformer overloading. The objective function for the OPF-based method is to minimise the number of EV disconnections. However, 1) there is no discussion of the computational complexity of each method in the paper, and 2) OPF-based method is not designed as a multiperiod form and the presented model does not have a storage model. Reference \cite{mehta_smart_2018} proposed two smart charging strategies with objectives of: 1) minimisation of total daily cost, and 2) peak-to-average ratio. The proposed strategies are tested with a 37-bus distribution system. The solution method is based on a heuristic-based method which is not fast enough for online operational strategies. Reference \cite{zhang_fast_2019} introduced a fast solving method for the coordinated charging of EVs based on linearisation of branch power flow. They implemented their work on the IEEE 33-bus distribution network, which is not a large case study. Reference \cite{shi_model_2019} developed a model predictive control-based approach to solve the joint problem of EV charging scheduling and power control. The objective is to minimise both EV charging cost and energy generation cost while satisfying the daily household and EV power demand. The largest test-case studied in the paper is the IEEE 118 bus. The proposed method is, however, not fast enough to handle large case studies. The authors of \cite{kotsalos_horizon_2019} proposed a horizon optimisation control framework in order to schedule the operation of the distribution network efficiently. The main objective of the proposed optimisation is to abide by the operational constraints to ensure secure operational scheduling. The operational constraints are voltage bounds and rated power bounds. The optimisation problem is based on multi-period three-phase Optimal Power Flow (OPF) which can be solved by a classic Non-Linear Programming (NLP) solver. However, the implemented method is tested in a small-scale Distribution Network (DN). The computational complexity of the proposed method is not discussed in the paper.\\
{Zaferanlouei et al. \cite{zaferanlouei_battpower_2021} proposed a tailored algorithm which is computationally efficient for the integration of EV into power distribution grid. The proposed method is applied on large-scale DN with the large population of EVs. This work is the next paper in a row with a more focus on deep analyses and discussions, which comprehends the impact of EV on distribution grid and how to mitigate it.}     

\subsection{Contributions and Paper Structure}
From a power system point of view, optimal charging of EVs is performed using a multiperiod ACOPF formulation in order to maximise the utilisation of renewable generation and cost-effective generation. In addition, market prices, system losses, as well as grid constraints such as line and transformer congestion, voltage levels and voltage angles are taken into consideration.\\
Although substantial efforts have been undertaken in order to propose optimised charging of EVs in the distribution grid, as reviewed in subsection \ref{EVCharge}, none of them suggests a tool that is computationally fast and demonstrated for a large-scale distribution grid.\\

The test case considered in this paper uses a large-scale distribution grid with maximum details available in Norway for analysis. To our knowledge, the work presented in this paper is the first ever attempt to do a comprehensive analysis of the impact of EV charging demand on real distribution grid.\\

Thus, the main novelty of this paper consists of considering optimal scheduling of EV charging while also including the full ACOPF formulation alongside operational constraints of the distribution grid in a real, large scale case. The main contribution of this work is the demonstration of the tractability and scalability of the model introduced in \cite{zaferanlouei_battpower_2021} in a real large-scale Norwegian distribution grid with 856 consumers, 974 buses, and 1023 lines. Following, we demonstrate how the proposed multiperiod ACOPF formulation can ensure stable grid operation by smartly charging EVs in a real distribution grid in Norway.\\

In section \ref{problemFormulate}, we briefly describe the MPOPF formulation and inputs of the BATTPOWER solver. Next, the input optimisation data are introduced and discussed in \ref{inputData}. The methodology of this study is elaborated in section \ref{methodology}. The numerical results are presented and further discussed in section \ref{results}. Finally we end the paper with concluding remarks.

%% file: tex/problemFormulate.tex
\section{Problem Formulations} \label{problemFormulate}
For a given power network with $n_b \in \mathbb{N}$ number of buses, $n_g\in \mathbb{N}$ number of generators, $n_l\in \mathbb{N}$ number of lines, $n_y\in \mathbb{N}$ number of storage devices (EV and ESS), a general MPOPF formulation with optimisation horizon of $T$ and time steps of $t=\{1,\cdots ,T\}$ can be written as: 
\begin{subequations}\label{multiperiod}
\begin{alignat}{4} 
&\min_{\mathbf{X}} F(\mathbf{X})\\ \label{multiperiodG.a}
\textrm{s.t. } G(\mathbf{X})&= \begin{bmatrix}
   \widetilde{G}(\mathbf{X}) \
   \overline{G}(\mathbf{X}) \
   \overline{G}^s(\mathbf{X})
\end{bmatrix}^\top&&=0  \in   \mathbb{R}^{N_{g} \times 1}\\\label{multiperiodH.b}
H(\mathbf{X})&=\begin{bmatrix}
   \widetilde{H}(\mathbf{X}) \quad
   \overline{H}(\mathbf{X})
\end{bmatrix}^\top&&\leq 0  \in   \mathbb{R}^{N_{h} \times 1} 
\end{alignat}
\end{subequations}
where the vector of total variables in the MPOPF problem $\mathbf{X} \in \ \mathbb{R}^{N_{x} \times 1}$  where $N_{x} = TN_{x_t}$, is shown in (\ref{X}) 
\begin{equation}
\begin{multlined}
\label{X}
\mathbf{X}= \big[\mathbf{x}_{1} \quad  \mathbf{x}_{2}\quad ...\quad \mathbf{x}_{t}\quad ...\quad  \mathbf{x}_{T}\big]^\top 
\end{multlined} 
\end{equation}
and the corresponding variables $\mathbf{x}_{t}$ for each time $t$ are defined through: 
\begin{equation}
\begin{multlined}
\label{eqn:x_t}
{\mathbf{x}_{t}= \big[\boldsymbol{\Theta}_{t} \  \boldsymbol{\mathcal{V}}_{t} \ \boldsymbol{\mathcal{P}}^{\mathrm{g}}_{t} \ \boldsymbol{\mathcal{Q}}^{\mathrm{g}}_{t} \   \boldsymbol{\mathcal{SOC}}_{t}\  \boldsymbol{\mathcal{P}}_{t}^\mathrm{ch} \ \boldsymbol{\mathcal{P}}_{t}^\mathrm{dch} \ \boldsymbol{\mathcal{Q}}_{t}^\mathrm{s}\big]^\top
 \ \mkern-10mu}_{1 \times N_{x_t}}
\end{multlined} 
\end{equation}
In this paper, besides  MATPOWER's well-known input matrices $\mathbf{BUS}$, $\mathbf{BRANCH}$, $\mathbf{GEN}$, and $\mathbf{GENCOST}$, we introduce new input matrices $\mathbf{BATT}$, $\mathbf{AVBP}$, $\mathbf{CONCH}$, $\mathbf{CONDI}$, $\mathbf{AVBQ}$, $\mathbf{AVG}$, $\mathbf{SOCi}$, $\mathbf{SOCMi}$ in order to capture the dynamic behaviour of MPOPF and especially large-scale integration of EVs. The description and purpose of each input are elaborated in \ref{appendixC} of this paper. In the following subsections, we expand the objective function and constraints of MPOPF.
\subsection{Objective Function}
The objective function of the entire optimisation period is the integral of objective functions for each period $\forall t \in \{1,\dots,T\}$. 
\begin{subequations}\label{eqn:OF}
\begin{alignat}{4}
&F(\mathbf{X})=f_{t=1}(\mathbf{x}_{1})+f_{t=2}(\mathbf{x}_{2})+\dots+f_{t=T}(\mathbf{x}_{T})\\
&f_{t}= (\boldsymbol{\mu}_t^{\mathrm{spot}})^\top\boldsymbol{\mathcal{P}}_t^{\mathrm{g}} 
\end{alignat}
\end{subequations}
where $ \boldsymbol{\mu}_t \in \mathbb{R}^{n_{b} \times 1}$ is the vector of the marginal price at time $t$. It contains the generators' marginal costs for the generators' buses and zero for the rest. $\boldsymbol{\mathcal{P}}_t^{\mathrm{g}}$ is the vector of active injected power by generator which is part of the variables in \eqref{eqn:x_t} at time $t$. 
\subsection{Constraints}
In this subsection, we extend the constraints brought in \eqref{multiperiod}. They are: a) $\widetilde{\mathbf{G}}(\mathbf{X})\in \mathbb{R}^{N_{gn} \times 1}$ is the vector of non-linear equality constraints, corresponding to the nodal power balance, AC power flow equations, and is extended in \eqref{19.a}.  b) $\overline{\mathbf{G}}(\mathbf{X})\in \mathbb{R}^{N_{gl} \times 1}$ is the vector of linear equality constraints, except the linear equality of storage devices, and is shown in \eqref{19.b}. c) $\overline{\mathbf{G}}^s(\mathbf{X})\in \mathbb{R}^{N_{gs} \times 1}$ is the linear equality of storage devices, and is shown in \eqref{19.c}. d) $\widetilde{\mathbf{H}}(\mathbf{X})\in \mathbb{R}^{N_{hn} \times 1}$ is the vector of nonlinear inequality constraints corresponding to line flows, and is clarified in \eqref{19.d}. Lastly, e) $\overline{\mathbf{H}}(\mathbf{X})\in \mathbb{R}^{N_{hl} \times 1}$ is the vector of linear inequality constraints known as box constraints, and is extended in \eqref{19.e}.
\begin{subequations}\label{19}
\begin{alignat}{4}
   \widetilde{\mathbf{G}}(\mathbf{X})=&
\begin{bmatrix}
    \widetilde{\mathbf{g}}(\mathbf{x}_{1}) \
    \widetilde{\mathbf{g}}(\mathbf{x}_{2}) \
    \dots \
    \widetilde{\mathbf{g}}(\mathbf{x}_{T}) 
\end{bmatrix}^\top \label{19.a}\\
   \overline{\mathbf{G}}(\mathbf{X})=&
\begin{bmatrix}
    \overline{\mathbf{g}}(\mathbf{x}_{1}) \
    \overline{\mathbf{g}}(\mathbf{x}_{2}) \
    \dots \
    \overline{\mathbf{g}}(\mathbf{x}_{T}) \
\end{bmatrix}^\top \label{19.b}\\
  \overline{\mathbf{G}}^s(\mathbf{X})=&
\begin{bmatrix}
    \overline{\mathbf{g}}^s(\boldsymbol{\tau}_1) \
    \overline{\mathbf{g}}^s(\boldsymbol{\tau}_2) \
    \dots \
    \overline{\mathbf{g}}^s(\boldsymbol{\tau}_T) 
\end{bmatrix}^\top \label{19.c} \\
   \widetilde{\mathbf{H}}(\mathbf{X})=&
\begin{bmatrix}
    \widetilde{\mathbf{h}}(\mathbf{x}_{1}) \
    \widetilde{\mathbf{h}}(\mathbf{x}_{2}) \
    \dots \
    \widetilde{\mathbf{h}}(\mathbf{x}_{T}) 
\end{bmatrix}^\top \label{19.d}\\
  \overline{\mathbf{H}}(\mathbf{X})=& 
\begin{bmatrix}
    \overline{\mathbf{h}}(\mathbf{x}_{1}) \
    \overline{\mathbf{h}}(\mathbf{x}_{2}) \
    \dots \
    \overline{\mathbf{h}}(\mathbf{x}_{T}) 
\end{bmatrix}^\top \label{19.e} 
\end{alignat}
\end{subequations}
where $N_{g} =N_{gn}+N_{gl}+N_{gs}$,\  $N_{gn} = Tn_{gn},\ N_{gl}= n_{{gl}_{t=1}}+n_{{gl}_{t=2}}+...+n_{{gl}_{t=T}},\ N_{gs}=Tn_y,\ N_{h} =N_{hn}+N_{hl},\ N_{hn}= Tn_{hn}, \ N_{hl}=n_{{hl}_{t=1}}+n_{{hl}_{t=2}}+...+n_{{hl}_{t=T}}+T(8n_{y})$, $\boldsymbol{\tau}_1 =\{\mathbf{x}_1\}$, $\boldsymbol{\tau}_t =\{\mathbf{x}_{t-1}$, $\mathbf{x}_t\}$ and $\mathcal{T}=\{\boldsymbol{\tau}_1$, $\boldsymbol{\tau}_2$, $\dots$, $\boldsymbol{\tau}_T\}= \{\{\mathbf{x}_1\}$, $\{\mathbf{x}_1$, $\mathbf{x}_2\}$,$\dots$,$\{\mathbf{x}_{T-1}$, $\mathbf{x}_{T}\}\}=\{\mathbf{x}_1$,$\mathbf{x}_2$,$\dots,\mathbf{x}_T\}$, thus  $\overline{\mathbf{G}}^s(\mathcal{T})$\footnote{Note that $\boldsymbol{\tau}_t$ and $\mathcal{T}$ are representations of two sets such that $\boldsymbol{\tau}_t\in \mathcal{T}$}$= \overline{\mathbf{G}}^s(\mathbf{X})$.
\subsubsection{Balance Constraint, (Full ACOPF)}
\begin{align}
\begin{split}
\label{active_Power_flow}
&\widetilde{\mathbf{g}}(\mathbf{x}_{t})\\
=&\begin{bmatrix}
\mathbf{C}_t^\mathrm{g} \boldsymbol{\mathcal{P}}_t^{\mathrm{g}}-\boldsymbol{\mathcal{P}}_t^{\mathrm{d}}-\mathbf{C}_t^\mathrm{ch} \boldsymbol{\mathcal{P}}_t^{\mathrm{ch}}+\mathbf{C}_t^\mathrm{dch} \boldsymbol{\mathcal{P}}_t^{\mathrm{dch}}-\Re{[\mathbf{\underline{S}}_t^{\mathrm{bus}}]} \\
\mathbf{C}_t^\mathrm{g} \boldsymbol{\mathcal{Q}}_t^{\mathrm{g}}-\boldsymbol{\mathcal{Q}}_t^{\mathrm{d}}+\mathbf{C}_t^\mathrm{s} \boldsymbol{\mathcal{Q}}_t^{\mathrm{s}}-\Im{[\mathbf{\underline{S}}_t^{\mathrm{bus}}]}
\end{bmatrix}=0
\end{split}
\end{align} 
\subsubsection{linear equality constraint \texorpdfstring{$\overline{\mathbf{g}}(\mathbf{x}_t)$}{}}
$\overline{\mathbf{g}}(\mathbf{x}_t)$ includes \eqref{eqn:equalitylineargrida}-\eqref{eqn:equalitylineargena} plus any other upper and lower bounds of variable $\mathbf{x}_t$ such that $x_t^\mathrm{min}=x_t^\mathrm{max}$, which can be user defined, and as such can be removed from the list of box constraints in \eqref{19.e} and is introduced here as a new linear equality \eqref{eqn:equalitylineargridd}.
\begin{subequations} \label{eqn:equalitylineargrid}
\begin{alignat}{6}
&\theta_t^{\mathrm{slack}}=0&&&\ \ \ \ &\label{eqn:equalitylineargrida}\\
 &p_{i,t}^\mathrm{ch}=0,\ \ \text{if } \quad \{\mathbf{AVBP}_{i,t}	\lor \mathbf{CONCH}_{i,t} \}=0\\
 & p_{i,t}^\mathrm{dch}=0,\ \text{if }\quad \{\mathbf{AVBP}_{i,t}	\lor \mathbf{CONDI}_{i,t} \}=0\\
&q_{i,t}^\mathrm{s}=0,\ \ \text{if }\quad \{\mathbf{AVBP}_{i,t}	\lor \mathbf{AVBQ}_{i,t}\}=0 \label{eqn:equalitylineargridc}\\
&p_{i,t}^\mathrm{g}=0,\ \ \text{if }\quad \mathbf{AVG}_{i,t}=0 \label{eqn:equalitylineargena}\\
&q_{i,t}^\mathrm{g}=0,\ \ \text{if }\quad \mathbf{AVG}_{i,t}=0 \label{eqn:equalitylineargenb}\\
 & x_t=x_t^\mathrm{min}= x_t^\mathrm{max} \  \ \text{if }  x_t^\mathrm{min}= x_t^\mathrm{max} \label{eqn:equalitylineargridd} 
 \end{alignat}
 \end{subequations}
 \subsubsection{Storage Device Constraints \texorpdfstring{$\overline{\mathbf{g}}^\mathrm{s}(\mathbf{x}_t)$}{}}
The vector of linear equality constraints corresponding to the storage devices $\overline{\mathbf{g}}^\mathrm{s}(\mathbf{x}_t)\in \ \mathbb{R}^{n_y \times 1}$ is defined from \eqref{storages}. 
\begin{equation}\label{storages}
 \mathbf{\overline{g}}^s(\boldsymbol{\tau}_t)=\mathbf{E}_t- \mathbf{E}_{t-1}-\mathbf{\Psi}^\mathrm{ch}\boldsymbol{\mathcal{P}}_t^\mathrm{ch}\Delta t+\frac {\boldsymbol{\mathcal{P}}_t^\mathrm{dch}\Delta t}{\mathbf{\Psi}^\mathrm{dch}} =0    
\end{equation}
where $\boldsymbol{\mathcal{SOC}}_t=\frac{\mathbf{E}_t}{\mathbf{E}^{max}}$ and $\{\boldsymbol{\mathcal{SOC}}_t$, $\boldsymbol{\mathcal{P}}_t^\mathrm{ch}$, $\boldsymbol{\mathcal{P}}_t^\mathrm{dch}$, $\boldsymbol{\mathcal{Q}}_t^\mathrm{s}\} \in \ \mathbb{R}^{n_y \times 1}$. Note that initial state of charge of each storage $i$ at time $t$ is defined as $e_{i,t-1}= e_{i}^{max}\mathbf{SOCi}_{i,t}$ where $\mathbf{SOCi}$ is the input matrix introduced in \ref{appendixC} such that an initial value of $\mathbf{SOCi}_{i,t}$ is allocated if one of the EV arrival conditions is satisfied: 1) $\mathbf{AVBP}_{i,t=1}=1$. 2) $\mathbf{AVBP}_{i,t-1}=0$ and $\mathbf{AVBP}_{i,t}=1$. 
\subsubsection{Line/Transformer Flow Constraints \texorpdfstring{$\widetilde{\mathbf{h}}(\mathbf{x}_{t})$}{}}
 Limitations on line/transformer flow, which are part of operational constraints, are now discussed. $\widetilde{\mathbf{h}}(\mathbf{x}_{t})$ is the vector of non-linear inequality constraints for time $t$ in Eq. \eqref{eqn:h_t}. 
 \begin{align}
\begin{split}
\label{eqn:h_t}
     \mathbf{\widetilde{h}}(\mathbf{x}_t)=
     \big[(\mathbf{\underline{S}}_t^{\mathrm{Line}})^*\mathbf{\underline{S}}_t^{\mathrm{Line}}-(\abs{\mathbf{\underline{S}}^{\mathrm{Line}}_\mathrm{max}})^{2}\big]
 \leq 0  \in \mathbb{R}^{n_{hn} \times 1}
 \end{split}
\end{align}
  \subsubsection{Box Constraints \texorpdfstring{$\overline{\mathbf{h}}(\mathbf{x}_{t})$}{}}
 Finally $\overline{\mathbf{h}}(\mathbf{x}_{t})$ is the set of box constraints of all the variables in \ref{X}. $N_{x_t} = n_{x}+4n_y$. Subscript $t$ stands for a specific time step in this paper. 
\begin{subequations}
\begin{alignat}{2}
&\boldsymbol{\Theta}^\mathrm{min}\leq \boldsymbol{\Theta}_t \leq \boldsymbol{\Theta}^{\mathrm{max}}\label{eqn:box1}\\ 
&\boldsymbol{\mathcal{V}}^{\mathrm{min}}\leq\boldsymbol{\mathcal{V}}_t\leq\boldsymbol{\mathcal{V}}^{\mathrm{max}} \label{eqn:box2}\\ 
&(\boldsymbol{\mathcal{P}}^{\mathrm{g}})^{\mathrm{min}}\leq\boldsymbol{\mathcal{P}}_t^{\mathrm{g}}\leq(\boldsymbol{\mathcal{P}}^{\mathrm{g}})^{\mathrm{max}}\label{eqn:box3}\\ 
& (\boldsymbol{\mathcal{Q}}^{\mathrm{g}})^{\mathrm{min}}\leq\boldsymbol{\mathcal{Q}}_t^{\mathrm{g}}\leq(\boldsymbol{\mathcal{Q}}^{\mathrm{g}})^{\mathrm{max}}\label{eqn:box4}\\
&\mathbf{SOCMi}_t\leq \boldsymbol{\mathcal{SOC}}_t \leq \boldsymbol{\mathcal{SOC}}^\mathrm{max}\label{eqn:boxe} \\ 
&(\boldsymbol{\mathcal{P}}^\mathrm{ch})^\mathrm{min} \leq \boldsymbol{\mathcal{P}}_t^{ch} \leq(\boldsymbol{\mathcal{P}}^\mathrm{ch})^\mathrm{max}\label{eqn:box6}\\ 
&(\boldsymbol{\mathcal{P}}^\mathrm{dch})^\mathrm{min} \leq \boldsymbol{\mathcal{P}}_t^\mathrm{dch}\leq (\boldsymbol{\mathcal{P}}^{\mathrm{dch}})^\mathrm{max}\label{eqn:boxg}\\
&(\boldsymbol{\mathcal{Q}}^\mathrm{s})^\mathrm{min} \leq
\boldsymbol{\mathcal{Q}}_t^\mathrm{s} \leq(\boldsymbol{\mathcal{Q}}^\mathrm{s})^\mathrm{max}\label{eqn:boxh} 
\end{alignat}
\end{subequations}
$\mathbf{SOCMi}_t$ is the vector of minimum state of charge taken from $\mathbf{SOCMi}$ input matrix for each time $t$. Inequality \eqref{eqn:boxe} is for the control of EV's state of charge before departure.

%% file: tex/inpData.tex
\section{Input Data} \label{inputData} 
The input data for the charge scheduling problem and BATTPOWER solver are described and discussed in this section. These data can be classified into the following categories:
\begin{enumerate}[(i.)]
\item Grid data
The test-case selected in this study is based on real case corresponding to a Norwegian distribution grid, illustrated in Fig. \ref{Norwegian1}. The test-case is originally presented and studied in \cite{zaferanlouei_integration_2017}. It has 32 medium-voltage (MV) 22kV to low-voltage (LV) 230V transformers. The entire system is fed from two buses: 1) a local generator on the left side of Fig. \ref{Norwegian1}, which is directly connected to the distribution grid through a 4kV/22kV transformer, and 2) the main feeder on the right side of it (bus no. 945), known as Point of Common Coupling (PCC). PCC is the connection to 66kV, and slack bus in this study. The test case consists of 974 buses, 1023 lines, and 856 consumers with hourly consumption data.\\
The grid data are translated into $\mathbf{BUS}$, $\mathbf{BRANCH}$, $\mathbf{GEN}$ and $\mathbf{GENCOST}$ matrices. The general format of these matrices are described in \ref{appendixA} and \ref{appendixC}.
\\
\afterpage{
\begin{figure*}[!htbp]
\centering
\includegraphics[width=4.5 in , height=3.2 in]{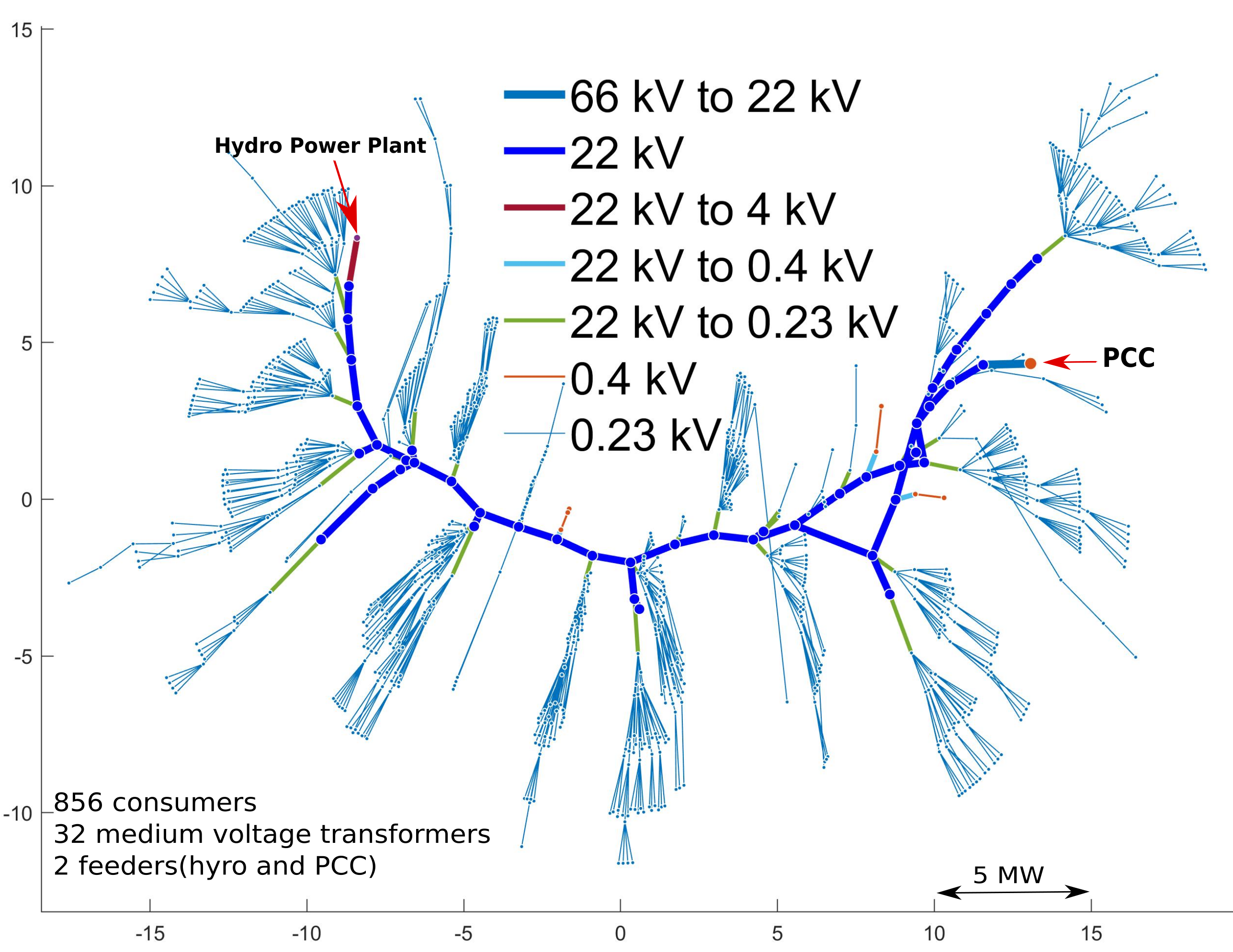}
\caption[]{Local distribution grid located in Norway with 856 costumers. Power transfer distance $PT_{ij}=\sum_{(i,j)\in \mathbf{BRANCH}} \left|PF_{ij}\right|$ visualisation technique is adopted from \cite{cuffe_visualizing_2017} to sketch this figure. The $PT_{ij}$ distance indicates how much aggregation of line utilisation occurs by the power transaction between each bus and the others\protect\footnotemark.}
\label{Norwegian1}
\end{figure*} 
\footnotetext{$PF_{ij}$ is the active power between the nodes i and j. For more information please see \cite{cuffe_visualizing_2017}}
}
\item Consumer's base-data: hourly base-load of 856 consumers are extracted through an algorithm described in Alg.\ref{Alg1}, subsection \ref{estimate}, located in \ref{appendixD} of this paper. 
\item EV data: EVs charge profile, time for  arrival and departure are elaborated in detail in subsection \ref{charge} located in \ref{appendixD} of this paper.
\end{enumerate}

%% file: tex/methodology.tex
\section{Methodology} \label{methodology}
In this section the methodology used is described. The general perspective towards the types of analysis conducted in this paper is given in the following subsections. The basis of proposed and applied control charge strategies are also explained in the next subsection. Two comprehensive cross-referenced tables (Table \ref{tab:Method} and Table \ref{tab:controlStrategy}) are presented to illustrate the concepts. 
\subsection{Overall perspective}
Three types of analysis are performed and the overview is shown in Table \ref{tab:Method}:
\begin{enumerate}[(A)]
\item \label{typeA} \textbf{DN Bottleneck:}\\
As described in the section \ref{intro}, the EV share will be substantially increased in the near future. Therefore, a comprehensive analysis of DN in terms of possible challenges are crucial. The first analysis goal is to simulate the current practice of EV charge strategy (Uncoordinated/dumb charge), while the EV share increases, and identify bottlenecks of DN in terms of congestion and voltage deviation. Power flow analysis is used in this analysis. The details of analysis, such as 15 min resolution of base-load of consumers (\ref{estimate}) and charging profile data (\ref{charge}) (i. EV drive distance, ii. EV distribution and energy consumption, iii. arrival and departure time and v. EV charge profile) are estimated according to Norwegian standards and reports (see \ref{appendixD} for more details). An overview of the simulation set-up can be found in Table \ref{tab:Method}. To investigate DN bottlenecks, we used the day of the year with peak base-load (12:00 PM 2 Feb. 2012 to 12:00 PM 3 Feb. 2012).
\item \label{typeB} \textbf{Price-Incentivised Charging:}\\
In brief, the analysis goal is to assess the impact of high price volatility on the price-incentivised charge strategy without considering operational constraints of the grid (MPOPF without network limits). It should be kept in mind that Norwegian energy price is quite flat\footnote{the main reason is that Norwegian power portfolio is mainly dominated by very flexible hydropower generation} (see Table \ref{tab:EV}). With the increasing share of varible renewable electricity production, the energy price is expected to become more volatile, as in Denmark with its high shares of wind turbines. Therefore, incorporating of volatile Danish energy price into the control charge strategy gives a perspective how price-incentivised charge is applicable for a near future scenario with higher share of wind and solar PV.
\item \label{typeC} \textbf{Price-Incentivised Charging with Network Limit:}\\
The purpose of analysis is to apply MPOPF with consideration of network constriants in a real Norwegian DN. We are able to make a socio-economical comparison/analysis with i. price-incentivised charge strategy (MPOPF without grid operational limits) and ii. uncoordinated/dumb charge strategy. A date with the highest peak power price in the year 2012, is selected for both coordinated and uncoordinated methods in order to make the comparison more accurate and insightful. Details regarding the analysis can be seen in Table \ref{tab:Method}.
\end{enumerate}

\begin{table}
\begin{center}
\begin{threeparttable}
\caption{Details about the type and purpose of analysis in this study.}
\label{tab:Method}
\small
\begin{tabular}{p{3cm}p{3.5cm}p{3.5cm}p{3.5cm}}
\toprule
Type of Analysis&  \ref{typeA}:  DN Bottleneck &   \ref{typeB}: Price-incentivised& \ref{typeC}: Price-incentivised with grid operational limit \\
\hline
\midrule
 Charge control strategy used in the analysis & Uncoordinated charge & Coordinated charge only MPOPF without network limits& \begin{tabular}{p{3.5cm}}i. Uncoordinated strategy\\ ii. Coordianted strategy both algortithms of MPOPF  with and without network limits\end{tabular}\\
\hline
Purpose of analysis &  {To investigate the DN bottlenecks (congestion and voltage), analysis result Table \ref{tab:dumbcharge}} & To study the impact of price volatiliy on the coordinated charge control strategy &  To compare the socio-economical benefits of \ref{typeA} and \ref{typeB}, results in Table \ref{tab:termination} \\
\hline
 Consumers'base-load (demand) &  24 hours with the   highest local DN power demand (MWh/h)  in the year 2012   &  24 hours with the   peak  power price (NOK/MWh) in the year 2012  &Similar to type \ref{typeB}\\
\hline
  Time resolution &  15 min (96 time period for 24 hours)   &  15 min (96 time period for 24 hours) & 15 min (96 time period for 24 hours)\\
\hline
  Date of data\tnote{1}\ in each analysis  & From 12:00 PM 2 Feb. 2012  to 12:00 PM 3 Feb. 2012  &  \begin{tabular}{p{3.5cm}} i. Base-load from noon to noon: 12:00 PM 1 Feb. 2012  to 12:00 PM 2 Feb. 2012\\ ii. Price signal of Danish simulation case  is selected form area price  of DK2  and one of the days with high price variations (high wind generation)\end{tabular}  &  \begin{tabular}{p{3.5cm}} i. Base-load from noon to noon: 12:00 PM 1 Feb. 2012  to 12:00 PM 2 Feb. 2012 \\ ii. Price signal from Trondheim price area  the same location   of real local DN  and the same date  from 12:00 PM 1 Feb. 2012  to 12:00 PM 2 Feb. 2012\end{tabular}  \\
\bottomrule
\end{tabular}
\begin{tablenotes}
\item[1] {\scriptsize The input data, which are specifically dependent on a date in this study, are i. consumer base-load and ii. price signal from power market.}
\end{tablenotes}
\end{threeparttable}
\end{center}
\end{table}
\subsection{Charge Scheduling Strategy}
The base case for the analysis is created with no EV presence to evaluate the impact caused by EV charging demand. A simple power flow analysis is conducted with the household's real demand. Next, system costs, system losses and daily aggregated power consumption are recorded for the further comparison. Three scenarios which are analysed and compared against the base case are:\\
\begin{enumerate}[(i.)]
\item \textbf{Uncoordinated/Dumb charge}\\
 The first scenario implements incremental population of EVs with uncoordinated (Dump) charging (the EVs will start to charge as and when they arrive). With the power flow analysis and comparison with the base case, the threshold population of EVs that causes problems for network operations are identified. Further, grid congestion and nodal voltage deviations/violations for the maximum possible EV populations for a Norwegian scenario, which is assumed to be 1.3 per household\footnote{for more detail about this assumption, please see \ref{vehiclenumber}, \ref{assumption} of this paper.}, are calculated. 
\item \textbf{Coordinated Charge Scheduling Strategy, MPOPF Without Network Limits}\\
The second scenario implements market price based optimised smart EV charging and includes power flow to calculate losses when it clears the nodal marginal prices. For the price based smart charging, two price regions are used: One for Norway where the test system is located and one from Denmark to assess more volatile prices. For more datails, see Table \ref{tab:controlStrategy}.
\item \textbf{Coordinated Charge Scheduling Strategy, MPOPF With Operational Limits}\\
The third scenario considers both market price and network constraints and does smart charging by running MPOPF with the objective of minimise the total operation cost of the modelled system. The thresholds of EV populations that cause congestion and voltage problem are identified for the third scenario. For more details, see Table \ref{tab:controlStrategy}.
\end{enumerate}
\subsection{BATTPOWER Application}
\begin{figure}[!htbp]
\centering
\includegraphics[width=3.5 in , height=1.8 in]{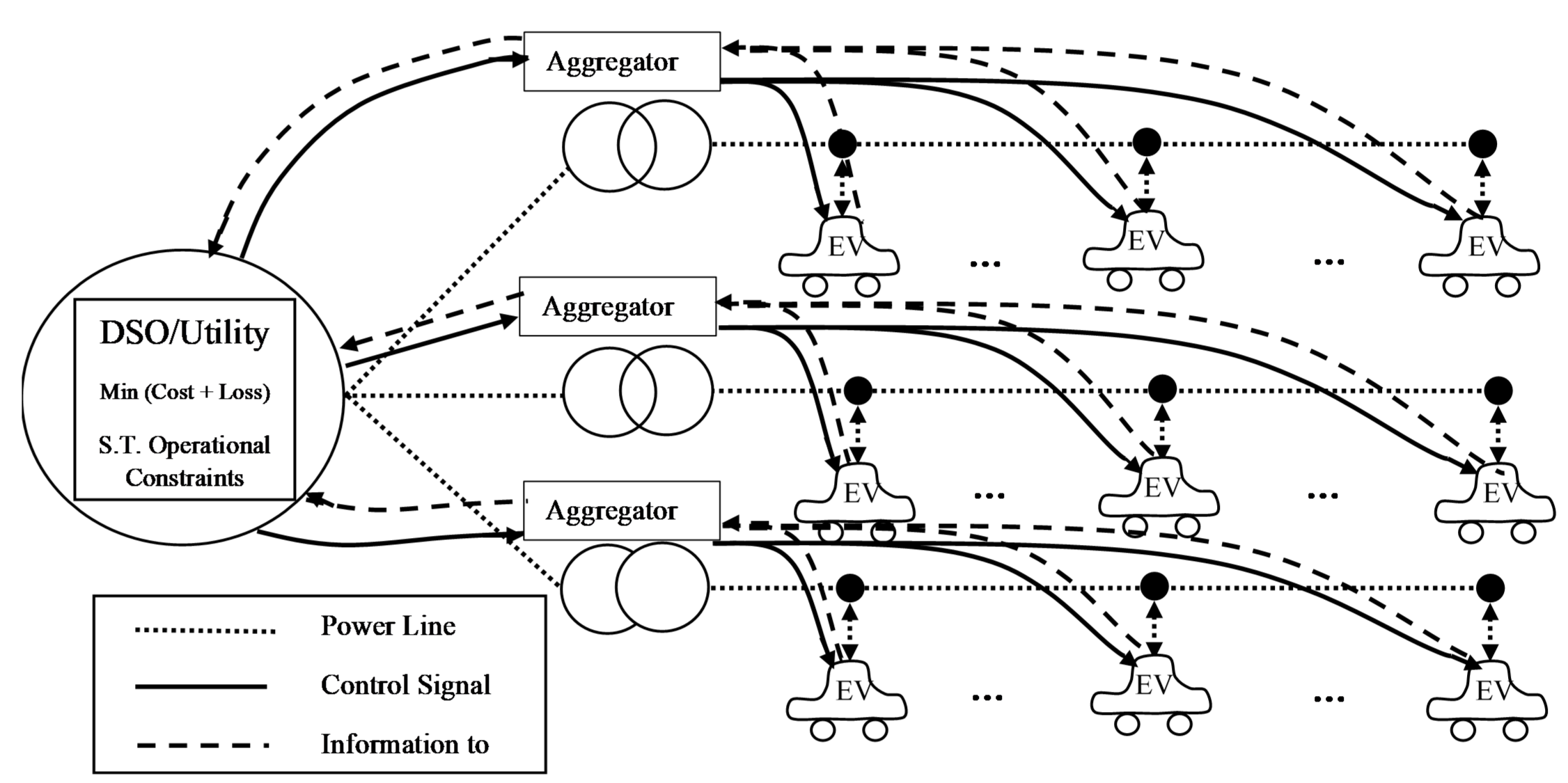}
\caption{Architecture of the proposed control charge scheduling strategies.} 
\label{Architecture} 
\end{figure}

The BATTPOWER solver developed in \cite{zaferanlouei_battpower_2021} is used in this study to simulate and incorporate the concept of coordinated EV charge scheduling algorithm. All input data (grid, price, base-load and EV) are fed into the BATTPOWER solver according to the details elaborated in \ref{inputData}, \ref{appendixC} and \ref{appendixD}. The solver is used to simulate two coordinated/optimal charge scheduling strategies.\\ 
The proposed control architecture for centralised charge scheduling of EVs is shown in Fig. \ref{Architecture}. The DSO(or charging operator) would follow three steps before running charging optimisation algorithm: 1) aggregate information regarding the arrival, departure, initial SOC, desired departure SOC and capacity of each EV's battery from the EV user, 2) predict the future base-load, and 3) predict the local generation of DERs (if the system adopts any). Next, the centralised optimisation algorithm is run and finally the EV charge demand scheduling would be cleared.\\
\begin{table}
\caption{Control strategies and their proposed simulation methods in this paper.}
\label{tab:controlStrategy}
\begin{center}
\small
\begin{tabular}{p{3cm}p{3cm}p{3cm}p{3cm}p{3cm}}
  \toprule
Control Strategy & \multicolumn{2}{c}{Uncoordinated (Dumb) charge} & \multicolumn{2}{c}{Coordinated charge}\\
\cmidrule(lr){1-1}
\cmidrule(lr){2-3}
\cmidrule(lr){4-5}
\cmidrule(lr){1-1}
\cmidrule(lr){2-3}
\cmidrule(lr){4-5}
Definition& \multicolumn{2}{c}{Charge on arrival}& \multicolumn{2}{c}{\begin{tabular}{p{5cm}}Charge occurs based on a control command form an algorithm\end{tabular}}\\
\cmidrule(lr){1-1}
\cmidrule(lr){2-3}
\cmidrule(lr){4-5}
\cmidrule(lr){1-1}
\cmidrule(lr){2-3}
\cmidrule(lr){4-5}
When can it be applied? & \multicolumn{2}{c}{Today's application}&Today's application &Near future\\
\cmidrule(lr){1-1}
\cmidrule(lr){2-3}
\cmidrule(lr){4-5}
\cmidrule(lr){1-1}
\cmidrule(lr){2-3}
\cmidrule(lr){4-5}
Proposed simulation model of the control strategy in this paper & \multicolumn{2}{c}{power flow} &  MPOPF without network constraints  &  MPOPF with network constraints\\
\cmidrule(lr){1-1}
\cmidrule(lr){2-3}
\cmidrule(lr){4-5}
\cmidrule(lr){1-1}
\cmidrule(lr){2-3}
\cmidrule(lr){4-5}
Objective function & \multicolumn{2}{c}{---} & Minimise accumulation of power price over optimisation horizon, see Eq.\ref{eqn:OF} also means: loss minimise  & Minimise accumulation of power price over optimisation horizon, see Eq.\ref{eqn:OF} also means: loss minimise \\
\cmidrule(lr){1-1}
\cmidrule(lr){2-3}
\cmidrule(lr){4-5}
\cmidrule(lr){1-1}
\cmidrule(lr){2-3}
\cmidrule(lr){4-5}
Constraints & \multicolumn{2}{c}{---} &\begin{tabular}{p{3cm}}i. Linear and nonlinear equality constraints: \eqref{19.a}-\eqref{19.c}\\ ii. Box constraints: \eqref{eqn:box3}-\eqref{eqn:boxh} \end{tabular}  & \begin{tabular}{p{3cm}}i. Power flow Eqs.\eqref{19.a}-\eqref{19.c}\\ii. Line/cable congestion \eqref{19.d}\\iii. Transformer congestion\eqref{19.d}\\iv. Voltage deviation \eqref{19.e} \\v. All box constraints \eqref{eqn:box1}-\eqref{eqn:boxh} \end{tabular}\\
\cmidrule(lr){1-1}
\cmidrule(lr){2-3}
\cmidrule(lr){4-5}
\cmidrule(lr){1-1}
\cmidrule(lr){2-3}
\cmidrule(lr){4-5}
Base-laod for the results &Highest consumers' load demand from 12:00 PM 2 Feb. 2012 to 12:00 PM 3 Feb. 2012& Highest peak price from 12:00 PM 1 Feb. 2012 to 12:00 PM 2 Feb. 2012& \multicolumn{2}{c}{\begin{tabular}{p{5cm}}Highest peak price from 12:00 PM 1 Feb. 2012 to 12:00 PM 2 Feb. 2012\end{tabular}}\\
\cmidrule(lr){1-1}
\cmidrule(lr){2-3}
\cmidrule(lr){4-5}
\cmidrule(lr){1-1}
\cmidrule(lr){2-3}
\cmidrule(lr){4-5}
 Results&\begin{tabular}{p{2cm}} DN bottleneck Table \ref{tab:dumbcharge}\end{tabular} & \begin{tabular}{p{3cm}} Result of uncoordinated charge simulation with different input data is used for a comparison with the proposed coordinated charging strategies: Table \ref{tab:termination} and \ref{tab:danishanlysis}\end{tabular}&  \multicolumn{2}{c}{\begin{tabular}{p{6cm}}Compared with base-load analysis in terms of: total energy consumption, loss, costs, savings, congestion managments,\\see: Table \ref{tab:termination} and \ref{tab:danishanlysis}\end{tabular}}\\
\bottomrule
\end{tabular}
\end{center}
\end{table}
\subsection{Assumption}\label{assumption}
All optimisations inputs, taken in this work, are derived and estimated based on some assumptions made, which can here be summarised:
\begin{enumerate}[(i.)]
\item Perfect foresight is used for MPOPF model (deterministic approach). 
\item Arrival, departure and daily travel distance are estimated based on methods elaborated in subsections \ref{charge}, \ref{distance}, \ref{consump}, \ref{arridep}, and \ref{profile}.
\item The hourly base-load data of the 856 registered consumers are estimated based on consumers' yearly energy consumption, i.e. the hourly load on two feeders for 8208 hours (32832 data for 15 min resolution), refer to subsection \ref{estimate}.
\item \label{zeroEV} The hourly consumer base-load, (extracted and estimated based on Alg.\ref{Alg1} in \ref{estimate}) are assumed to be no EV base-load (0\% EV penetration base-load). It should be noted that the original data of base-load are belong to the year 2012 with 0.4\% EV penetration; therefore, the assumption, made in this study, is very valid\footnote{Total number of registered EVs in 2012 are 9565 \cite{passengercars2012}, total number of registered passenger cars in 2012 are 2442960 \cite{EVs2012}. Thus, by the end of 2012, the penetration of EV in Norway is about 0.4 \%}. 
\item The efficiency of the EV charging inverter is constant.
\item A Real-Time Pricing (RTP) scheme is used, and the hourly energy tariffs are set equal to the Elspot day-ahead prices in Trondheim for 2012.
\item The average number of vehicles owned by a consumer is set to 1.3. This assumption is made due to the average number of vehicles per capita \cite{noauthor_statens_nodate} and the average number of people per household in mid-Norway \cite{noauthor_tabell_nodate}. Thus, the case study contains in total $856 \times 1.3=1113$ EVs.\label{vehiclenumber}
\item Cost functions of PCC and generator are similar and are a linear function of $f(\boldsymbol{\mathcal{P}}^{\mathrm{g}^{\mathrm{PCC}}},\boldsymbol{\mathcal{P}}^{\mathrm{g}^{\mathrm{gen}}})=\boldsymbol{\mu}^{\top}(\boldsymbol{\mathcal{P}}^{\mathrm{g}^{\mathrm{PCC}}}+ \boldsymbol{\mathcal{P}}^{\mathrm{g}^{\mathrm{gen}}})$ where $\boldsymbol{\mu} \in \mathbb{R}^{T\times 1}$ is the marginal hourly spot price (NOK/MW).  We assumed that the feeder and generator have similar hourly cost functions.
\item The EV batteries does not discharge power to the grid: $\boldsymbol{\mathcal{P}}_t^\mathrm{dch}=0$, in other words, minimum and maximum bounds are selected to be zero in the optimisation for all EVs and entire time horizon: $(\boldsymbol{\mathcal{P}}^\mathrm{dch})^\mathrm{min} = (\boldsymbol{\mathcal{P}}^\mathrm{dch})^\mathrm{max}=0$ 
\item The EV reactive power provision is zero: 
$\boldsymbol{\mathcal{Q}}_t^\mathrm{s}=0$, the same as item above: $(\boldsymbol{\mathcal{Q}}^\mathrm{s})^\mathrm{min}= (\boldsymbol{\mathcal{Q}}^\mathrm{s})^\mathrm{min}=0$
\item Initial SOC of EVs are introduced with input matrix of $\mathbf{SOCi}$ such that an initial value of $\mathbf{SOCi}_{i,t}$ is allocated if one of the arrival conditions are satisfied: 1) $\mathbf{AVBP}_{i,t=1}=1$. 2) $\mathbf{AVBP}_{i,t-1}=0$ and $\mathbf{AVBP}_{i,t}=1$.
\item Minimum SOC at time $t$ is defined through input matrix of $\mathbf{SOCMi}$, where constraint \eqref{eqn:boxe} is assumed to hold.
\end{enumerate}

%% file: tex/result.tex
\section{Results} \label{results} 
In this section, the numerical results of the proposed charge scheduling control of EVs are further presented. The input data are fed into the BATTPOWER solver through input matrices described in \ref{appendixC}. All assumptions made to conduct this study are listed in \ref{appendixD}, subsection \ref{assumption}. All cases shown in this section are performed on a workstation equipped with an Intel(R) Core(TM) i7-8650U CPU \@ 1.90GHz 16.GB RAM.\\
We simulate four cases to show how the MPOPF smart charge scheduling can be used for the effective utilisation of grid by avoiding grid congestions and voltage violations. Case 1 and 2 are with together simulation of DN to assess accommodation of a growing share of EVs. They are particularly designed to fulfill analysis type \ref{typeA}, assessment and analysis of the impact of uncoordinated control of charge scheduling of EVs. The case 3 focuses on price based active control and the case 4 focuses on MPOPF based active control.  The cases are listed below and then described in each subsection:
\subsection{Base-Load analysis}\label{base_load_analysis}
\begin{figure}[!htbp]
\centering
\includegraphics[width=3 in , height=4.5 in]{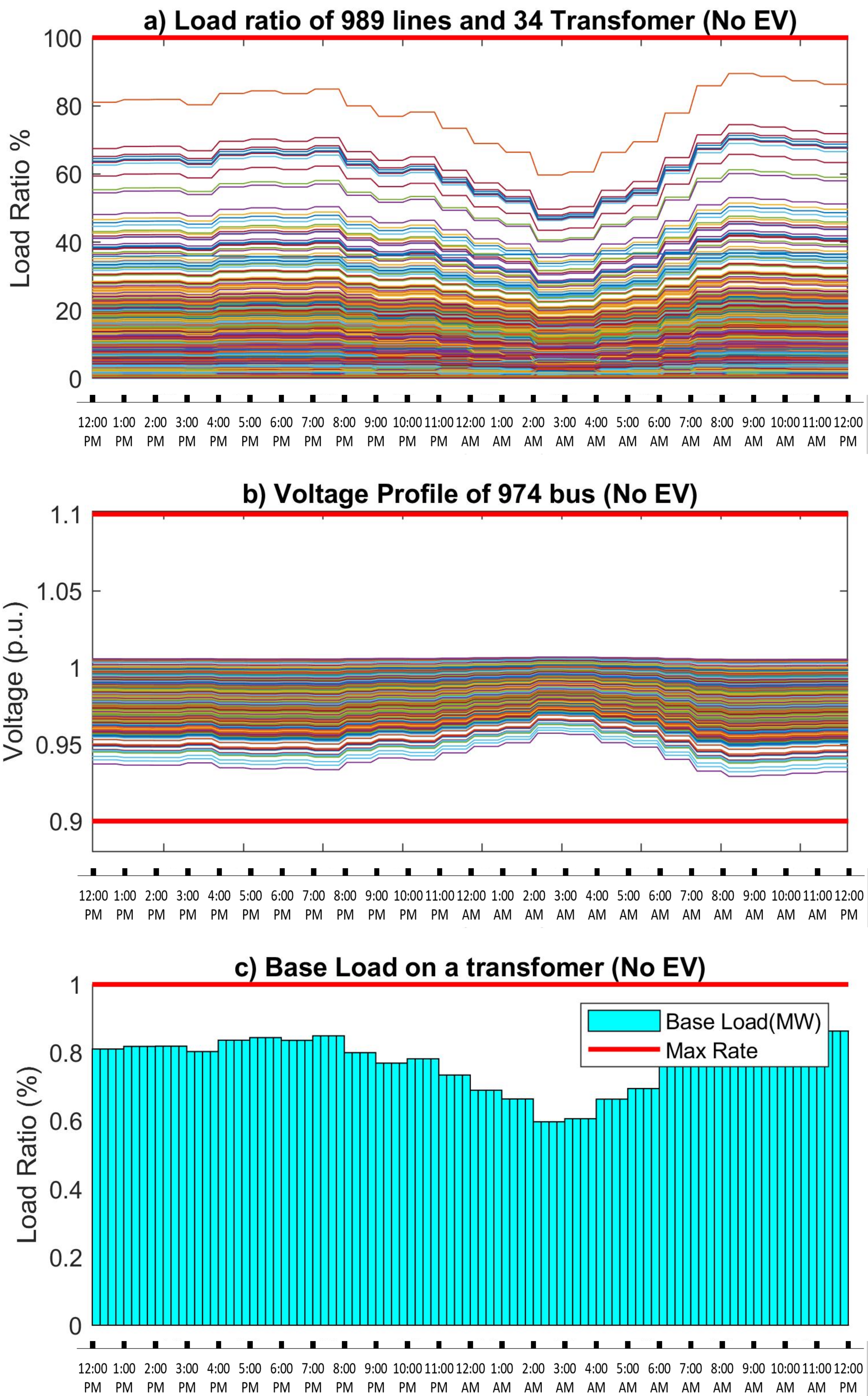}
\caption{power flow simulation of base-load demand. a) 15 min time resolution of load ratio of all lines/transformers, b) voltage profiles of 974 buses, outcome of load flow analysis, and c) Load ratio (\%) of the transformer no. 14 (22kV to 230 V) with 35 consumers.} 
\label{baseload}
\end{figure}

In order to analyse the impact of the base-load demand on the distribution grid, load flow analysis is conducted for the entire period of 342 days. No EV is inegrated to the grid, and the day with the highest load demand is selected (from 12:00 PM 2 Feb. 2012 until 12:00 PM 3 Feb. 2012). Fig. \ref{baseload} illustrates the results of base-load case, where: a) the load ratio of 989 lines and 34 transformers, b) the voltage profile of 974 buses, and c) the load ratio of the transformer  no. 14 (T54320) with 35 consumers, which is loaded close to its capacity in comparison with other transformers. There are two peaks which belong to the base-load demand of evening and morning in Fig. \ref{baseload} a) and c) respectively. No operational constraint of the grid is violated for the entire simulation period. This is the case of the highest base-load demand for Norwegian grid design.
\subsection{Uncoordinated Charging}
\begin{figure}[!htbp]
\centering
\includegraphics[width=3 in , height=4.5 in]{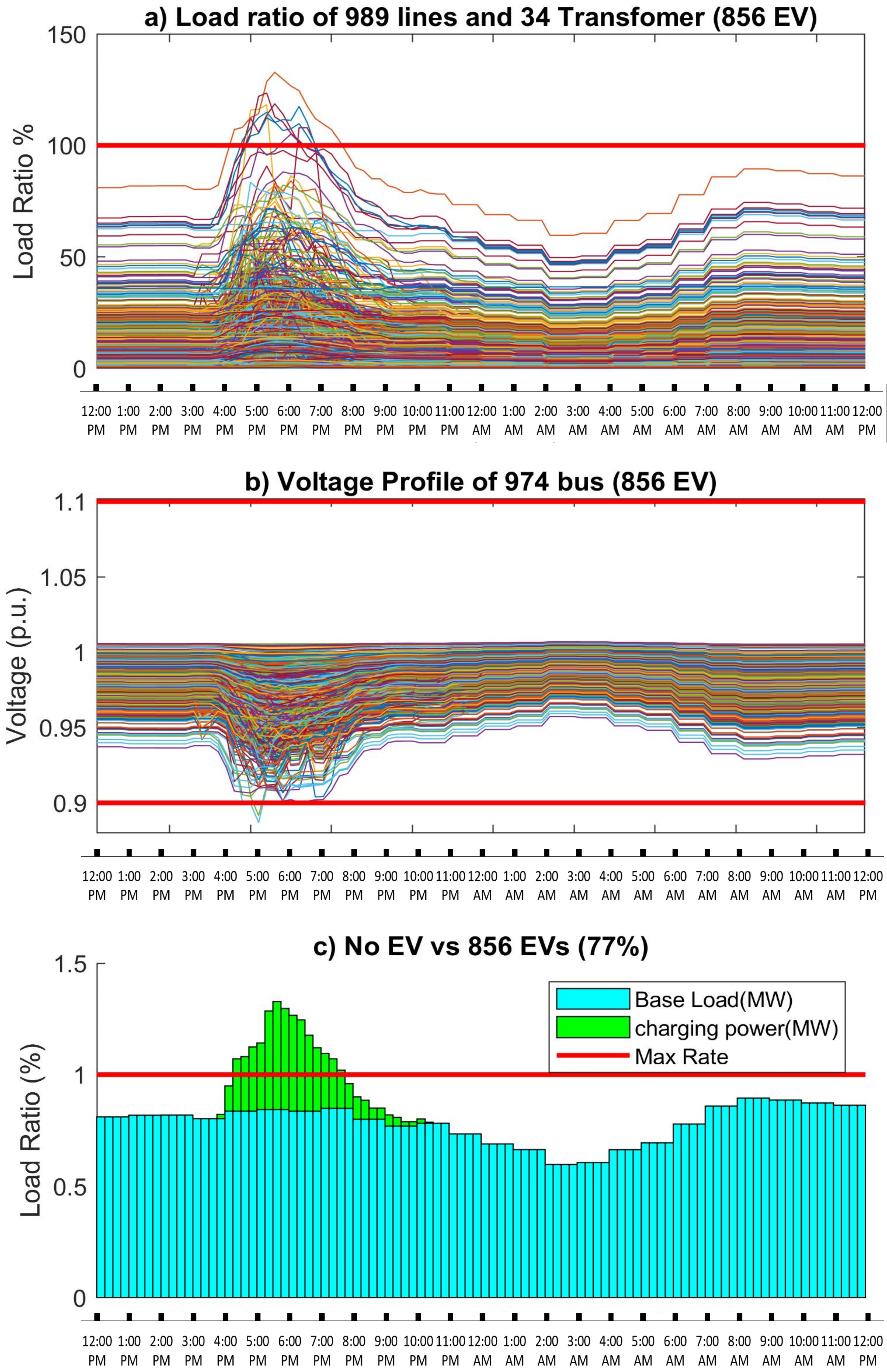}
\caption{power flow simulation of the uncoordinated 856 EVs charging in the local distribution grid. a) shows the load ratio (\%) of entire system components. b) voltage profile of all buses. c) Load ratio (\%) of the transformer no. 15 (22kV to 230 V) with only base-load and with uncoordinated EV's charge demand. The transformer is overloaded with a dumb-charging strategy.} 
\label{dumbCharge} 
\end{figure}

Uncoordinated (dumb) charging is referred to the charge scheduling control of EVs such that when an EV arrives, it connects to the grid until it is fully charged (SOC 100\%). Load flow analysis is conducted to simulate the impact of uncoordinated EV charge. The EV charge profile data are generated based on some assumptions elaborated in \ref{charge}.\\
Fig. \ref{dumbCharge} depicts the uncoordinated EV charge scheduling strategy for the case of 856 EVs in the mid-Norway distribution grid from 12:00 PM 2 Feb. 2012 until 12:00 PM 3 Feb. 2012. Fig. \ref{dumbCharge} a) shows 15-min resolution of load ratio of 989 lines and 34 transformers where load ratio on some lines/transformers are violated the maximum grid operational limits between the hours 17:00 and 19:00. b) illustrates the  15-min resolution of the voltage profile of 974 buses. Voltage fluctuation goes below the minimum voltage constraint of 0.9 p.u.. c) shows the load ratio of transformer no. 15.\\
{
\subsubsection{Analysis \ref{typeA}}
The simulations of uncoordinated charge scheduling of EVs are conducted for the cases of 100\footnote{The number of EVs in the DN} (9\%)\footnote{Penetration of EVs, as a percentage}, 200 (18\%), 300 (27\%), 500 (45\%), 600 (54\%), 856 (77\%), and 1113 (100\%) in the simulated area. Table \ref{tab:dumbcharge} summarises the load flow analysis results, for the uncoordinated EV charge scheduling strategy, where EV penetration increases from 0\% until 100\%. The capacity constraint reaches its maximum limit when the number of EVs is between 200-300 (18\%-27\%), but the voltage constraints are not violated. The voltage constraint is violated when the  number of EVs is in a range between 600-856 (54\%-77\%). The violation of line/transformer capacity limit is shown with double red lines and the violation voltage is shown with double blue lines in Table \ref{tab:dumbcharge}. Note that maximum voltage deviations are defined here as $1\ p.u.\pm 10 \%$.\\}
\begin{table}[htbp!]
\caption{Load flow analysis of the uncoordinated EV charge scheduling strategy. Bottlenecks of Norwegian DN. With uncoordinated charge of EV, the congestion happens when the share of EV is between 18\%-27\% (shown with double red lines), the voltage violation occurs when the share of EV is between 54\%-77\% (depicted with double blue lines).}
\begin{center}
\label{tab:dumbcharge}
\begin{threeparttable}
\begin{tabular}{l c c} 
\toprule
penetration & overloading (\%) \tnote{1}& voltage (p.u.) \tnote{2}\\
\midrule \midrule
Base-load & 89 \% & 0.93  \\
\hline
100 (9\%) & 94\%  &  0.927 \\  
\hline
200 (18\%)& 99\% & 0.921 \\
\arrayrulecolor{red}\hline\hline
300 (27\%) & 105 \%&  0.92 \\
\arrayrulecolor{black} \hline
500 (45\%) &  121 \%& 0.91 \\
\hline
600 (54\%) & 122\%& 0.91 \\
\arrayrulecolor{blue}\hline\hline
856 (77\%) & 140\%&  0.885\\
\arrayrulecolor{black}\hline
1113 (100\%) &150 \%& 0.875\\
\bottomrule
\end{tabular}
\begin{tablenotes}
\item[1] {maximum line/transformer overloading for 342 days}
\item[2] {minimum voltage in the simulated period, 342 days}
\end{tablenotes}
\end{threeparttable}
\end{center}
\end{table}
\subsection{Coordinated Charging without Considering the Operational Limits of the Grid} \label{without}
\begin{figure}[!htbp]
\centering
\includegraphics[width=3.5 in , height=6.2 in]{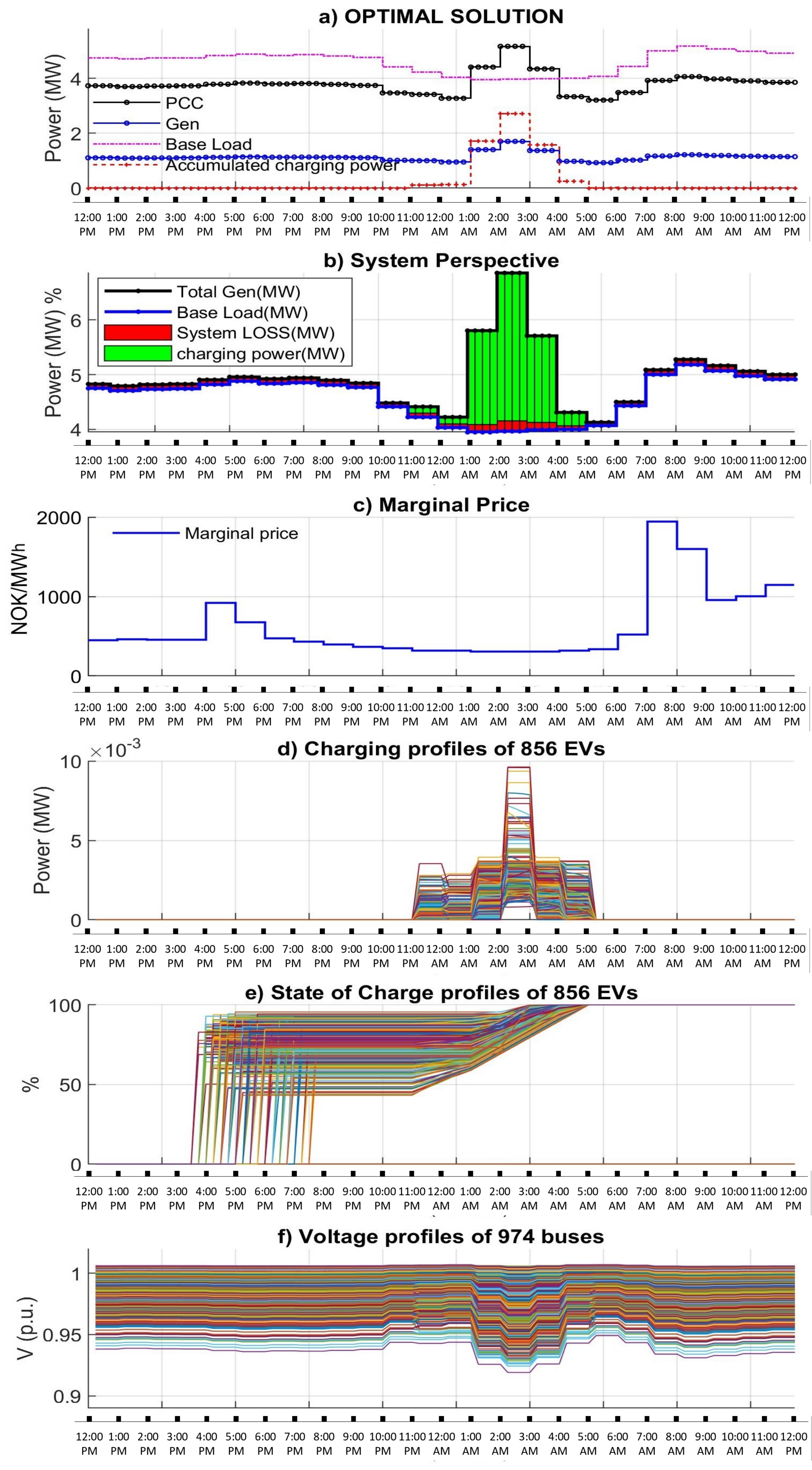}
\caption{The proposed centralised EV charge scheduling without grid operational constraints. 856 EVs charging in the local distribution grid. a) Aggregated load demand, production of generator and PCC import, b) residual load vs aggregated production. EV charge and loss can be seen here. c) daily profile of spot price d) charge profile of 856 EVs e) SOC of 856 EVs, the hard constraint to charge all of them to 100\% SOC. f) voltage profile of all buses.} 
\label{optimalCharge} 
\end{figure}

\begin{figure}[!htbp]
\centering
\includegraphics[width=3.5 in , height=6.2 in]{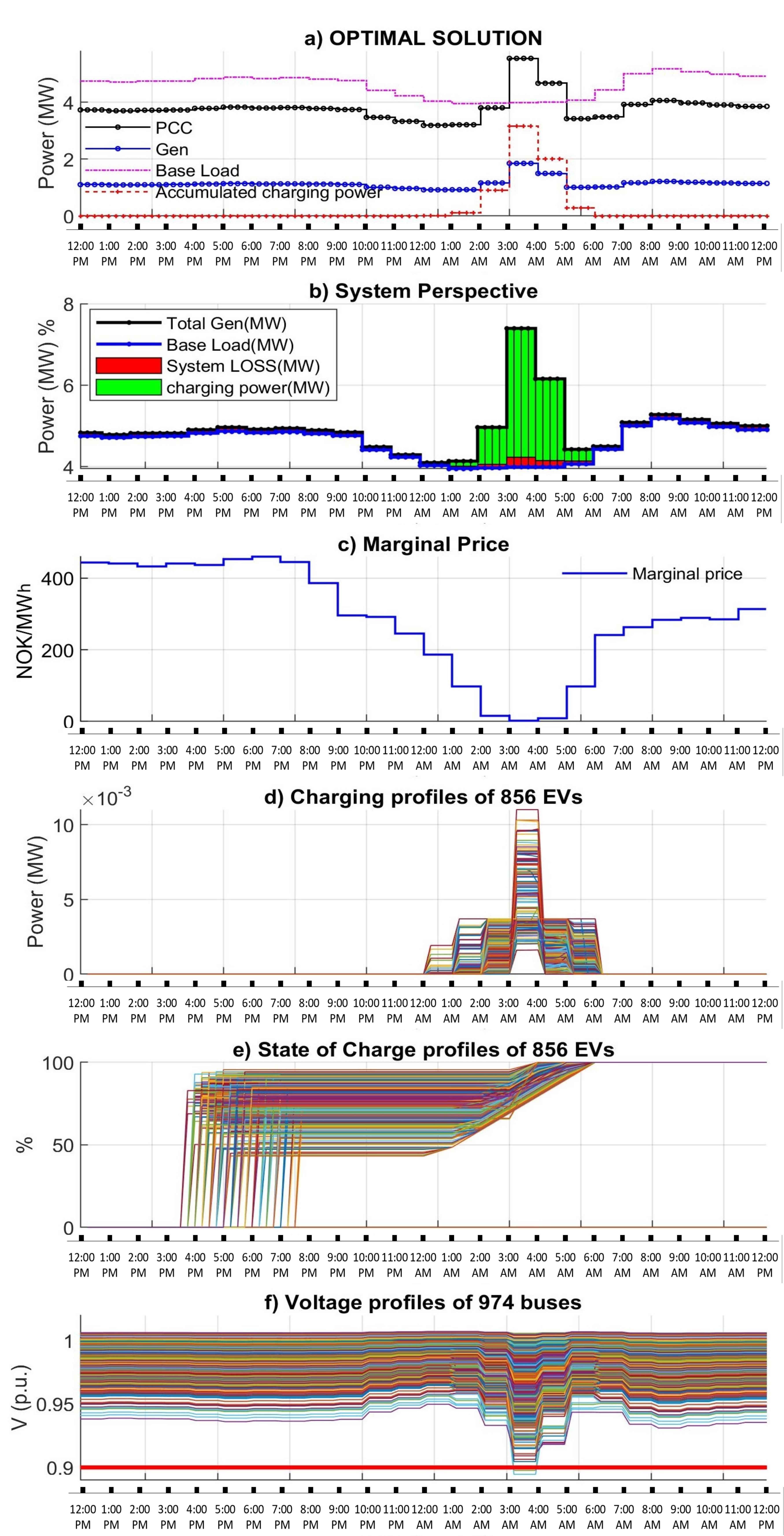}
\caption{Load-flow analysis of EV charge scheduling method. The base-load and EV penetration are similar to Fig. \ref{optimalCharge}. However, the spot price from DK2 (one of the days with the large price variations) is adopted to show that MPOPF without grid operational constraints could be volatile to inputs of the optimisation framework, since it forces the charge to occur during the lowest price time.} 
\label{Operational} 
\end{figure}

The proposed charge scheduling model, in this subsection, is a model which works based on a signal from DSO (or EV charge operator) to EV owners to react to: 1) the spot price, and 2) power losses in the distribution grid. The objective is to minimise the total cost by shifting the charging to low price intervals. The model respond to power losses, as  higher the power loss, the greater the total costs will be. No grid operational constraint is considered in this model. Fig. \ref{optimalCharge} shows the simulation results of the model. It is assumed that, the EV owner is equipped with a charger that react to the price signals. The control variables are considered to be the time of charge and rate of charge. In this subsection and the next subsection \ref{with}, the base-load and price input data are corresponding to a period between 12:00 PM Feb. 1 2012 and 12:00 PM Feb. 2 2012, in which the highest spot price - 2000 (NOK/MWh) - of the year 2012 occurs at 8:00 AM Feb. 2. Fig. \ref{optimalCharge} a) shows residual base-load, production of the generator and net energy import at PCC, by assuming that $\boldsymbol{\mu}$ (cost of energy) is the same for both (the one produce by the generator and the imported energy at PCC from the upstream network). b) shows the total generation versus total base-load. In between generation and consumption plots, loss and EV charge demand are shown in red and green coloured bar plots. c) shows spot price from \cite{nordpool}. d) illustrates the proposed EV charge scheduling demand of 856 EVs (77\%). e) shows the SOCs of 856 EVs, when they arrive with initial SOC of $\mathbf{SOCi} \in \mathbb{R}^{n_y \times T}$. The distribution of EV arrival can here be seen. f) depicts the voltage profile of 974 buses with a resolution of 15 min. There is no voltage violation, since the Norwegian spot price is flat during midnight. This is the opposite case in Fig. \ref{Operational} e) where the Danish spot price is adopted.\\
{
\subsubsection{Analysis \ref{typeB}}
The Norwegian price profile is often flat during midnight, as depicted in Fig. \ref{optimalCharge}, part (c) due to hydro dominated electricity generation. However, this is not the case for the Danish price profile as shown in Fig. \ref{Operational} part (c). Therefore, the same charge scheduling strategy (coordinated charge using MPOPF without considering network limits), is repeated with only one difference that the price signal is adopted from one of the days with large price variation of Danish price profile of DK2 as shown in Fig. \ref{Operational}, part (c).\\}

{
The impact of price volatility can be analysed with two different criteria:
\begin{enumerate}[I.]
\item \textbf{voltage violation:} The Danish case shows a voltage drop below 0.9 p.u. during the lowest price, between 3:00 AM and 4:00 AM, (see Fig. \ref{Operational}, part (f)) due to simultaneous charge scheduling of many EVs (rebound effect). The adoption of the Danish price profile is done for to assess a near future scenario with the growing share of intermittent renewable energy resources. In such a scenario, a stable price profile such as hydro dominated Norwegian price profile with small STD can be considered a rare case. The volatile Danish electriciy price area of DK1 and DK2 are two distinct examples of a near future scenario, where the price sometimes becomes negative. Therefore, the adoption of volatile DK2 price signal demonstrates that a charge scheudling control strategy based on locational marginal pricing is not an applicable/sufficient control charge strategy to charge a large population of EVs (as an example: in Fig. \ref{Operational}, part (c)). It should be kept in mind that, based on MPOPF without operational limits of the grid, only i. marginal price production of generators, and ii. DN loss reflect on locational marginal pricing, and not violation of operational constraints of the grid.  Note that the simulation of near future scenario, here in this analysis, is valid with the assumption of no reinforcement of DN.\\  
\item \textbf{congestion:} Line/transformer congestion happens in the both cases no matter of high or low price variations which can be seen in Fig. \ref{Optimal2} part (a).
\end{enumerate}
} 
\subsection{Coordinated Charging with Operational Limits of the Grid} \label{with}
\begin{figure}[!htbp]
\centering
\includegraphics[width=3.2 in , height=2.5 in]{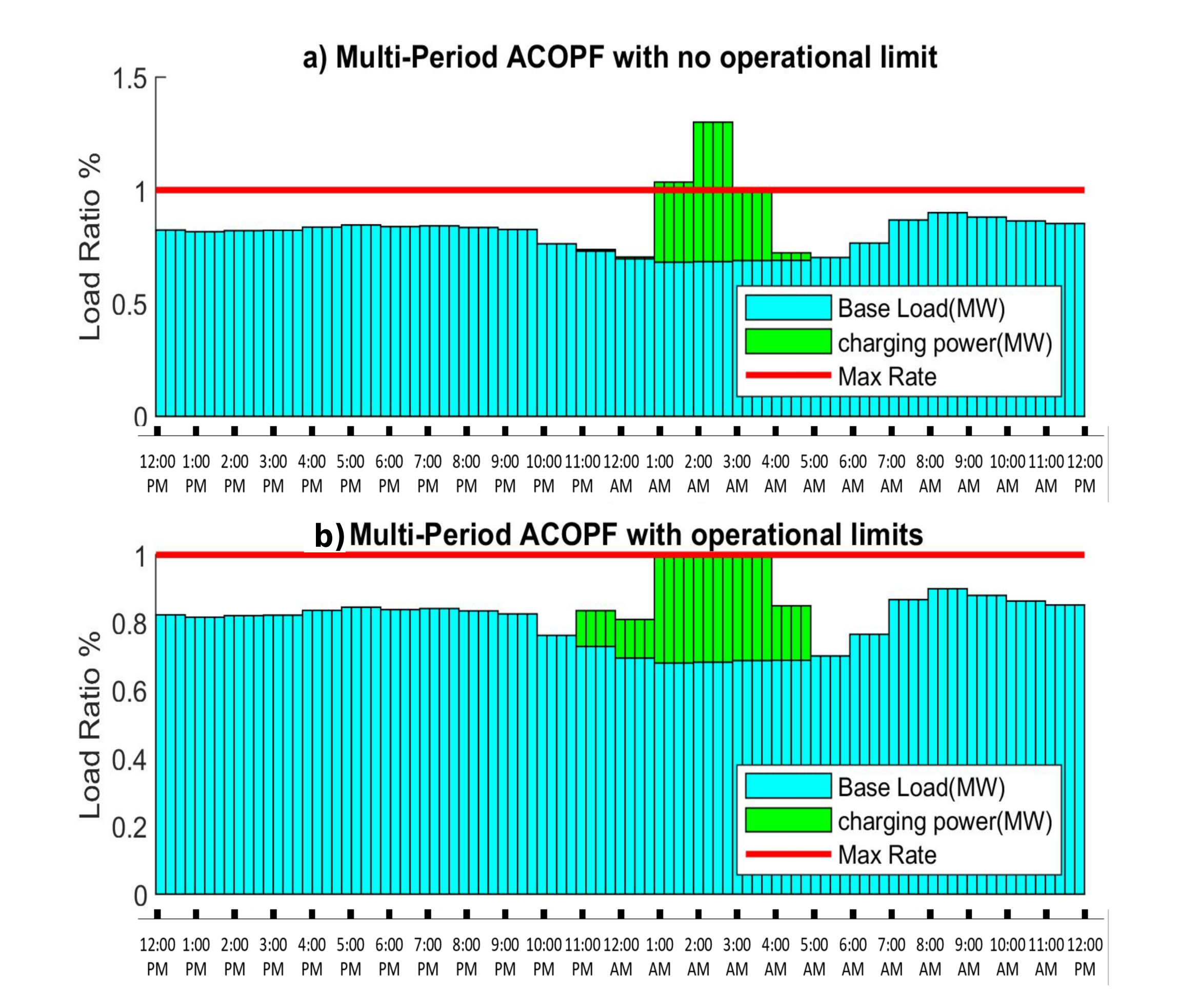}
\caption{The load ratio of the transfromer no. 14. The entire DN hosting 856 EV users in the simulation, where 35 of them are located and fed by this transformer. Comparison of: a) MPOPF without grid operational limit, and b) MPOPF with grid operational limit. Both cases are the load ratio of the transformer no. 14 (22kV to 230 V).} 
\label{Optimal2} 
\end{figure}

The centralised coordinated EV charge scheduling, considering operational constraints of line/transformer overloading, and voltage violation. The objective function is similar to previous subsection \ref{without} where the goal is to minimise the system cost plus system loss subject to not violating operational constraints of the grid. The difference between MPOPF without grid constraints \ref{without} and with grid constraints \ref{with} are depicted in Fig. \ref{Optimal2}. As can be seen, the EV charge demand profile is sharp in Fig. \ref{Optimal2} a). However, EV charge demand of profile shown in b) is flattened. The reason for this is the overloading constraint on the transformer. EV load is shifted to the next  time intervals.\\
{
\subsubsection{Analysis \ref{typeC}} \label{resultsanalysisc}
Table \ref{tab:termination} compares three different strategies of EV charge scheduling, presented in this section. These strategies are examined for: 1) daily energy consumption, 2) system loss, 3) system cost, 4) daily/yearly saving, and finally 5) robustness of each method to charge EVs without interruption. Daily energy consumption \footnote{This calculation is done with the assumption that the components overloading are ignored.} (MWh) and system loss are almost similar in the three proposed cases. However, system cost (NOK) is slightly different, as it is possible to save 2.6 \% with strategy 2 and 3 respectively.\\
It should be kept in mind that: 1) these values are only based on energy price (summation of hourly energy consumed multiplied to hourly spot price for the period under study), which means in reality this value might be much higher than the one presented here\footnote{Note that the net value of enegry bought from energy market, which is calculated in this study, is only about 34\% of consumer's bill \cite{eriksen_rme_nodate}}, and 2) the saving value, is calculated subject to the assumption that the spot price is flat during midnight in the Norwegian system, as can be seen \cite{nordpool}. If in the simulation results, shown in Fig. \ref{optimalCharge}, we adopt a 24-hour Danish price profile, the saving values shown in Table \ref{tab:termination} would have significantly increased. In order to prove this, the results shown in Table \ref{tab:termination} are re-simulated with the Danish area price of DK2\footnote{All the input data are similar with the only exception of the input signal price is taken from DK2, shown in Fig. \ref{Operational}, part (c)} and re-drawn in the Table \ref{tab:danishanlysis}, in \ref{appendixF}. As the results indicates in the table, the daily saving values are more than 8\%. This means with a volatile signal price, the coordinated charge sheduling strategies have more savings and thus the higher incentives to be implemented.\\
Lastly, as can be seen, lines/transformers are overloaded by the EV penetration around 20\% for dumb-charge strategy and 36\% for MPOPF without network constraints. However, the last MPOPF with grid operational limits extends the EV charging schedule until the departure time of EV. Thus it can handle (1113) 100\% EV penetration in the distribution grid, which is not possible with the other cases.\\}
\begin{table*}[htbp!]
\begin{center}
\caption{[a) total energy production, b) active system loss, and c) system cost] in three different operational mode.}
\label{tab:termination}
\begin{threeparttable}
\begin{tabular}{l c c c c c c} 
\toprule
 {\scriptsize Method} & \begin{tabular}[c]{@{}c@{}}{\scriptsize Daily\tnote{4} Energy}\\ {\scriptsize Consumption (MWh)}\end{tabular}  &  \begin{tabular}[c]{@{}c@{}}{\scriptsize Active Loss}\\ {\scriptsize (MWh)}\end{tabular}&  \begin{tabular}[c]{@{}c@{}}{\scriptsize System Cost}\\ {\scriptsize (NOK)}\end{tabular} &               \begin{tabular}[c]{@{}c@{}}{\scriptsize Daily Saving}\\ {\scriptsize (NOK)--(\%)}\end{tabular}  & \begin{tabular}[c]{@{}c@{}}{\scriptsize Yearly Saving}\\ {\scriptsize (NOK)}\end{tabular}  &  \begin{tabular}[c]{@{}c@{}}{\scriptsize Max EV}\\ {\scriptsize hosting Capacity}\end{tabular} \\  
\midrule \midrule
 \tnote{1}&118.83 & 2.24 & 75,927.1 & -- & -- &220 EV (20\%) \\
\midrule
 \tnote{2} & 118.74&2.15 & 73,973.5& 1,953.2-- 2.572 \%  & 712,916& 400 EV (36\%) \\
\midrule
 \tnote{3}&118.74&2.15 & 73,974.5 & 1,952.7 -- 2.571 \% & 712,743& 1113 EV (100\%) \\
\bottomrule
\end{tabular}
\begin{tablenotes}
{\footnotesize
\item[1]  Dumb Charging. 
\item[2]  MPOPF without grid operational limits. 
\item[3]  MPOPF with grid operational limits. 
\item[4]  Obtained based on the base-load input data of date of 12:00 PM 1 Feb. 2012 to 12:00 PM 2 Feb. 2012.}
\end{tablenotes}
\end{threeparttable}
\end{center}
\end{table*}
\section{Discussion}
{
\subsection{Impact of EVs on the distribution grid, as a rule of thumb}
Uncoordinated/dumb charge simulation of a real distribution grid is conducted in this paper, using real aggregated base-load data of 856 consumers. There are some small but important points one could take from the real data analysis parts of this paper.\\
Table \ref{tab:thumb} shows two calculated simple criteria in order to provide a general perspective for power systems researchers and engineers, the ratio of average EV load on average base-load and the ratio of the maximum of aggregated EV load on the average aggregated base-load. The former value is 3\% while the latter value is computed to be 24\%. The main question here is how to interpret these two values/criteria in a meaningful and simple manner. The first one (the average load of EV per average base-laod) means ratio of the EV load on the base-load, in general. In another word, 3\% means that the EV load is not much in compare with overall base load and can be ignored. However, the second criterion with the value of 24\% can be interpreted differntly. Although the general ratio of EV load on base-load is 3\%, the ratio of the maximum (maximum of value over time) aggregated EV load for the local community per the average (mean over time) aggregation of base-load for the same local community is considerable. In another word, EV's impact on the grid can occur for a short period of time.}
\begin{table}
\centering
\begin{threeparttable}
\caption{The impact of EVs on the distribution grid, as a rule of thumb.}
\label{tab:thumb}
\footnotesize
\begin{tabular}{p{2cm}p{2cm}p{2cm}p{2cm}p{2cm}p{1.5cm}}
\toprule
Data Type  &Dimension 1: Power (Load)  &Dimension 2: Time   &\begin{tabular}{c}Operation\end{tabular} &\begin{tabular}{c} Source\\ of Data\end{tabular} &Ratio \\
\hline
\midrule
EV Data ${\begin{bmatrix} p_{i,j}\end{bmatrix}\ \mkern-10mu}_{856 \times 35136}$ & $i \in$\ 856 EV load data  &$j \in\ $35136 timestamp (15 min resolution for 365 days)\tnote{*} & \begin{tabular}{p{2cm}}Average over both $\{i,j\} =$ \begin{tabular}{c}0.083\end{tabular}\end{tabular} &Simulated based on reports and Norwegian living standards, see: \ref{charge}\tnote{**}  & \multirow{12}{*}{\begin{tabular}{p{1.5cm}} $\frac{0.083}{2.8}$\\  = 0.03\\:= 3\% \end{tabular}}\\
\cmidrule(lr){1-5}
Base-Load ${\begin{bmatrix} l_{i,j}\end{bmatrix}\ \mkern-10mu}_{856 \times 32832}$ & $i \in $ \ active base-load of 856 consumers  & $j \in $\ 32832 timestamp (15 min resolution)\tnote{\#} & \begin{tabular}{p{2cm}}Average over both $\{i,j\}=$\begin{tabular}{c}2.8\end{tabular} \end{tabular} &Extracted from real aggregated data, see \ref{estimate}\tnote{\#\#}& \\
\hline
\hline
Aggregated EV Data ${\begin{bmatrix}\sum_{i=1}^{856} p_{i,j}\end{bmatrix}\ \mkern-10mu}_{1\times 35136}$  &\begin{tabular}{c} ---\end{tabular} &similar to *   &$\underset{j}{max}{\begin{bmatrix}\sum_{i=1}^{856} p_{i,j}\end{bmatrix}}$ $=568$kW &similar to ** &\multirow{6}{*}{\begin{tabular}{p{1.5cm}} $\frac{568}{2400}$\\  = 0.24\\:= 24\% \end{tabular}}\\
 \cmidrule(lr){1-5}  
 Aggregated Base-Load${\begin{bmatrix}\sum_{i=1}^{856} l_{i,j}\end{bmatrix}\ \mkern-10mu}_{1\times 32832}$&\begin{tabular}{c} ---\end{tabular}   & similar to\#  & $\underset{j}{mean}{\begin{bmatrix}\sum_{i=1}^{856} l_{i,j}\end{bmatrix}}$ $=2400$ kW&similar to \#\# &\\
\bottomrule
\end{tabular}
\begin{tablenotes}
\item[1] {\scriptsize Cross-reference symbols in the table cells are: *, **, \# and \#\#.}
\end{tablenotes}
\end{threeparttable}
\end{table}
\subsection{Coordinated charge and the need for a high computational toolbox}
An application of a high-performance solver is presented. The solution proposed is to solve full multiperiod ACOPF equations coupled with storage device coupling constraints represented as EVs all together as a centralised charge scheduling framework. \\
The computational time to solve the Norwegian distribution grid case study, with the specifications of 974 buses, 1023 lines, 2 generators, 856 consumers, 1113 EVs (100\%), and 24-hours optimal planning horizon (resolution of 96 time steps in the optimisation horizon) is 790 second, with the current formulations presented in this paper. The presented method could be potentially a handy tool used by DSO and for a near future scenario.\\
In the current market design, the end users buy electricity from retailer companines, where retailer participate in the spot market to buy their predicted demand power and sell it as an average price to the end users. A portion of consumers' monthly bill is allocated as the grid tariff. In the whole mechanism, the active transmission power loss is estimated and DSO is charged with a penalty for losses in the distribution grid. Therefore, with the current market design, DSO is not only responsible for ensuring the safe and stable operation of the distribution grid, but also technically responsible for the efficient operation of the distribution grid. In addition, charge scheduling of EV, is not only about DSO and its current defined operational tasks, but also economical incentives. The energy bought from the spot market in order to satisfy EV charge demand can be imported during low electricity price, as indicated in this study.\\
Although through coordinated EV charge, centralised MPOPF algorithm, cost minimisation could be an ultimate objective, the simulation results here suggest that through the implementation of a centralised algorithm, system cost is minimised. Thus, it should be kept in mind that the relative saving is not significant in comparison with a dumb-charge strategy (maximum 2.4\%). A main question is brought to mind in this respect: Who is going to use MPOPF tools like the one presented here in practice? Considering grid operation, only the DSO would be able to run an optimisation because it is the only stakeholder with access to the grid data. However, the DSO is not allowed to buy and trade electricity as it violates current regulation (It is not technically a market player).\\ 
Some research studies \cite{sortomme_optimal_2011, yao_hierarchical_2013, ortega-vazquez_electric_2013, chung_electric_2019} suggested that a new market player, called an aggregator (or EV operator), is required in order to control and schedule EV charge demand.  A more viable approach in a market design perspective is that the aggregator should act as the market player to coordinate the asset management.\\ 
The drawback is that the aggregator's task is to satisfy their customers which are the asset owners. In addition, the aggregator would normally have no access to grid data. We therefore get two possible market frameworks that support wide-scale use of the tool presented in this paper. The DSO runs MPOPF based on expected available flexibility and places bids in the local flexibility market to try to achieve the results from this tool. Alternatively, an aggregator controls the assets, but also has access to grid data. The advantages of the latter option could be discussed as the aggregator could be a suitable market player in this respect. It could receive information from EV owners and send command back to the EV, what rate to charge and when to charge. The proposed centralised EV charge scheduling could be a running algorithm in this market player. Minimising system cost and loss w.r.t. satisfying grid constraint could be a feasible option. In order to coordinate demand profile such that the rate of charge and time of charge can be controllable. In this respect, a contract between end user and DSO could be placed such that the aggregator could send the control signal to end user and control the charge scheduling process.     

%% file: tex/conclusion.tex
\section{Conclusion}
We presented an EV charge scheduling algorithm for large-scale integration  of EV in the distribution grid, minimising energy costs plus system loss subject to grid operational constraints to ensure a safe and reliable distribution grid operation. The proposed method incorporates multiperiod AC optimal power flow (MPOPF) coupled with energy storage device constraints. A large-scale distribution grid is chosen as a benchmark in this study, in order to assess different strategies which could potentially be implemented by DSO. Moreover, the presented EV charge scheduling method performs with a highly computational efficiency. We showed that a cost minimisation function adopts a larger numbers of EVs in the benchmark case study. However, neither the uncoordinated, nor cost minimisation strategy are capable of integrating 100\% EV in the local distribution grid. Therefore, the DSO has two options to fully electrify the transport sector; either to strengthen the distribution grid or apply a charge scheduling EV mechanism. The presented study suggests a fast, scalable, control of EV charge mechanism as a solution to a future sustainable electricity grid. To our knowledge, the work presented in this paper is the first ever attempt to do a comprehensive analysis of the impact of EV charging demand on real distribution grid. The inference of the analysis says that the Norwegian distribution networks are more prone to congestion problems than the voltage problems for the EV demand.

%% file: tex/acknow.tex
\section*{Acknowledgments}
This work has been supported by the project consortium of ``Grid and Charging Infrastructure of the Future  FuChar, grant number (295133/E20)'', ``Intelligent dispatching and optimal operation of cascade hydropower plants based on spatiotemporal big data-IntHydro, grant number: (309997/E20)'' and  ``CINELDI  Centre for intelligent electricity distribution, an 8-year Research Centre under the FME scheme (Centre for Environment-friendly Energy Research, grant number: 257626/E20)'' all funded by the Research Council of Norway.

%% file: tex/appA.tex
\section{Nomenclature} \label{appendixA}
\mbox{}
\setlength{\nomlabelwidth}{3cm}

\nomenclature[A, 01]{LV}{Low Voltage}
\nomenclature[A, 02]{MV}{Medium Voltage}
\nomenclature[A, 03]{LV}{High Voltage}
\nomenclature[A, 04]{MPOPF}{MultiPeriod AC Optimal Power Flow}
\nomenclature[A, 05]{DER}{Distributed Energy Resource}
\nomenclature[A, 06]{DSO}{Distribution System Operator}
\nomenclature[A, 07]{EV}{Electric Vehicle}
\nomenclature[A, 08]{PCC}{Point Common Coupling}
\nomenclature[A, 09]{SOC}{State of Charge}
\nomenclature[A, 10]{DN}{Distribution Network}
\nomenclature[A, 11]{AVE}{Average}
\nomenclature[A, 12]{STD}{Standard deviation}

\nomenclature[B, 01]{$f_{t}, F$}{Objective function of time step $t$, objective function of entire period refer to Eq. \eqref{eqn:OF}.}
\nomenclature[B, 02]{$\mathbf{g}, \mathbf{G}$}{Vector of equality constraint of one time step and vector of equality constraint of entire optimisation horizon $T$.}
\nomenclature[B, 03]{$\mathbf{h}, \mathbf{H}$}{Vector of inequality constraint of one time step and vector of inequality constraint entire optimisation horizon $T$.}
\nomenclature[B, 04]{$\mathbf{E}_t$}{Vector of energy of $n_y$ storage devices and EVs (MWh).}
\nomenclature[B, 05]{$\boldsymbol{\mu}_t^{\mathrm{spot}}$}{Vector of marginal hourly spot price (NOK/MWh).}
\nomenclature[B, 06]{$\mathbf{\underline{S}}^{\mathrm{Line}}_t$}{$\in \mathbb{C}^{2n_l \times 1}$ Vector of rated capacity (MVA) of $2n_l$ Lines (from and to) at time $t$.}%
\nomenclature[B, 07]{$\mathbb{B}$}{A binary matrix/set.}
\nomenclature[B, 08]{$\mathbb{C}$}{A complex matrix/set.}
 \nomenclature[C, 01]{{$\mathbf{BUS}$}}{{Matrix of $\mathbf{BUS}$ contains input data w.r.t. the  buses (in terms of types, loading, initial voltage magnitude and angle in OPF calcaulations, and etc) of the network.}}
 \nomenclature[C, 02]{{$\mathbf{BRANCH}$}}{{Matrix of $\mathbf{BRANCH}$ contains input data w.r.t. the lines (bus connections). These can be line resistance (p.u.), reactance (p.u.), charging susceptance (p.u.), line operational limits, and etc.}}
 \nomenclature[C, 03]{{$\mathbf{GEN}$}}{{Matrix of $\mathbf{GEN}$ contains specification of generators and their operational status.}}
 \nomenclature[C, 04]{{$\mathbf{GENCOST}$}}{{Matrix of $\mathbf{GENCOST}$ contains generation cost function.}}
 \nomenclature[C, 05]{{$\mathbf{BATT}$}}{{Marix of $\mathbf{BATT}$ contains battery specifiction and its location in the network.}}
 \nomenclature[C, 06]{{$\mathbf{AVBP}$}}{{Matrix of $\mathbf{AVBP}$ contains the connection status (availability) of the battery.}}
 \nomenclature[C, 07]{{$\mathbf{CONCH}$}}{{Matrix of $\mathbf{CONCH}$ contains the conditions for charging (0 not allowed, 1 allowed).}}
 \nomenclature[C, 08]{{$\mathbf{CONDI}$}}{{Matrix of $\mathbf{CONDI}$ contains the condition for discharging (0 not allowed, 1 allowed).}}
 \nomenclature[C, 09]{{$\mathbf{AVBQ}$}}{{Matrix of $\mathbf{AVBQ}$ contains the condition for battery reactive power control (0 not allowed, 1 allowed).}}
 \nomenclature[C, 10]{{$\mathbf{AVG}$}}{{Matrix of $\mathbf{AVG}$ contains the condition for generator power control (0 not allowed, 1 allowed).}}
 \nomenclature[C, 11]{{$\mathbf{SOCi}$}}{{Matrix of $\mathbf{SOCi}$ contians the intitial state of charge at different time.}}
 \nomenclature[C, 12]{{$\mathbf{SOCMi}$}}{{Matrix of $\mathbf{SOCMi}$ contains the minimum state of charge restriction value.}}
\nomenclature[C, 13]{{$\mathbf{PD}$}}{{Matrix of $\mathbf{PD}$ contains time series of active loads on the buses.}}
\nomenclature[C, 14]{{$\mathbf{QD}$}}{{Matrix of $\mathbf{QD}$ contains time series of reactive loads on the buses.}}
\nomenclature[C, 15]{{$\texttt{BATT\_BUS}$}}{{Column of $\texttt{BATT\_BUS}$ indicates the location of ESS/EV per bus number.}}
\nomenclature[C, 16]{{$\texttt{SOC\_OPT}$}}{{Column of $\texttt{SOC\_OPT}$ indicates the initial value of variable of$\boldsymbol{\mathcal{SOC}}_{t}$ in the optimisation algorithm.}}
\nomenclature[C, 17]{{$\texttt{PCH\_OPT}$}}{{Column of $\texttt{PCH\_OPT}$ indicates the initial value of variable of $\boldsymbol{\mathcal{P}}_{t}^\mathrm{ch}$ in the optimisation algorithm.}}

\nomenclature[C, 18]{{$\texttt{PDICH\_OPT}$}}{{Column of $\texttt{PDICH\_OPT}$ indicates the initial value of variable of $\boldsymbol{\mathcal{P}}_{t}^\mathrm{dch}$ in the optimisation algorithm.}}

\nomenclature[C, 19]{{$\texttt{Q\_INJ\_OPT}$}}{{Column of $\texttt{Q\_INJ\_OPT}$ indicates the initial value of variavle of $\boldsymbol{\mathcal{Q}}_{t}^\mathrm{s}$ in the optimisation algorithm.}}
\nomenclature[C, 20]{{$\texttt{MBASE}$}}{{Column of $\texttt{MBASE}$ indicates the total MVA base of this ESS/EV.}}

\nomenclature[C, 21]{{$\texttt{EFF\_CH}(\mathbf{\Psi}^\mathrm{ch})$}}{{Column of $\texttt{EFF\_CH}(\mathbf{\Psi}^\mathrm{ch})$ indicates the efficiency of charge of ESS/EV.}}
\nomenclature[C, 22]{{$\texttt{EFF\_DICH}(\mathbf{\Psi}^\mathrm{dch})$}}{{Column of $\texttt{EFF\_DICH}(\mathbf{\Psi}^\mathrm{dch})$ indicates the efficiency of discharge of ESS/EV.}}
%
%
%
%
%
%
%
%
%
%
%
%
\nomenclature[D, 01]{$n_b, n_g, n_l, n_y$}{Number of buses, generators, branches, and storage devices/EVs.}
\nomenclature[D, 02]{$N_x, N_{x_t}$}{Total number of variables, number of variables at time $t$.}
\nomenclature[D, 03]{$N_g, N_{gn}, N_{gl}, N_{gs}$}{Total number of equality constraints, number of nonlinear equality constraints (balance ACOPF constraints), number of linear equality constraints except storage devices and EVs, number of linear equality constraints of storage devices and EVs.}
\nomenclature[D, 04]{$N_{h}, N_{hn}, N_{hl}$}{Total number of inequality constraints, number of nonlinear inequality constraints (line flow), number of linear inequality constraints.}
\nomenclature[D, 05]{$n_{gn}, n_{gl_t}$}{Number of nonlinear equality constraints (balance ACOPF constraints) of time $t$, number of linear equality constraints of time $t$ (storage devices and EVs are not included).}
\nomenclature[D, 06]{$n_{hn}, n_{hl_t}$}{number of nonlinear inequality constraints (line flow) of time $t$, number of linear inequality constraints of time $t$.}
\nomenclature[D, 13]{$T$}{Number of steps in the optimisation horizon.}
\nomenclature[D, 14]{$\boldsymbol{\Theta}^\mathrm{min}, \boldsymbol{\Theta}^\mathrm{max}$}{Minimum and maximum limit of the voltage phase angles of $n_b$ buses, (rad).}
\nomenclature[D, 14]{$\boldsymbol{\mathcal{V}}^{\mathrm{min}}, \boldsymbol{\mathcal{V}}^{\mathrm{max}}$}{Minimum and maximum limit of the voltage magnitude of $n_b$ buses (p.u.).}
\nomenclature[D, 14]{$(\boldsymbol{\mathcal{P}}^{\mathrm{g}})^{\mathrm{min}}, (\boldsymbol{\mathcal{P}}^{\mathrm{g}})^{\mathrm{max}}$}{Minimum and maximum limit of the active power generated of $n_g$ generators (MW).}
\nomenclature[D, 14]{$(\boldsymbol{\mathcal{Q}}^{\mathrm{g}})^{\mathrm{min}}, (\boldsymbol{\mathcal{Q}}^{\mathrm{g}})^{\mathrm{max}}$}{Minimum and maximum limit of the reactive power generated of $n_g$ generators (MVA).}
\nomenclature[D, 14]{$\boldsymbol{\mathcal{SOC}}^{\mathrm{min}}, \boldsymbol{\mathcal{SOC}}^{\mathrm{max}}$}{Minimum and maximum limit of the state of charge of $n_y$ storage devices and EVs.}
\nomenclature[D, 15]{$(\boldsymbol{\mathcal{P}}^\mathrm{ch})^\mathrm{min}, (\boldsymbol{\mathcal{P}}^\mathrm{ch})^\mathrm{max}$}{Minimum and maximum limit of the rated charging capacity of $n_y$ storage devices and EVs (MW).}
\nomenclature[D, 15]{$(\boldsymbol{\mathcal{P}}^\mathrm{dch})^\mathrm{min}, (\boldsymbol{\mathcal{P}}^\mathrm{dch})^\mathrm{max}$}{Minimum and maximum limit of the rated discharging capacity of $n_y$ storage devices and EVs (MW).}
\nomenclature[D, 15]{$(\boldsymbol{\mathcal{Q}}^\mathrm{s})^\mathrm{min}, (\boldsymbol{\mathcal{Q}}^\mathrm{s})^\mathrm{max}$}{Minimum and maximum limit of the rated reactive power capacity of $n_y$ inverters of storage devices and EVs (MVA).}
\nomenclature[D, 17]{$\Delta t$}{time step.}%
\nomenclature[D, 19]{$\mathbf{\Psi}^\mathrm{ch}, \mathbf{\Psi}^\mathrm{dch}$}{Efficiency vectors of charge and discharge of $n_y$ number of storage devices and EVs which are taken from input matrix of $\mathbf{BATT}$.}
\nomenclature[D, 22]{$\mathbf{\underline{Y}}^{\mathrm{fr}}, \mathbf{\underline{Y}}^{\mathrm{to}}, \mathbf{\underline{Y}}^{\mathrm{Line}}$}{Line admittance matrix of from bus $i$ to bus $j$ of $n_l$ number of lines and $n_b$ number of buses, line admittance matrix of to bus $i$ from bus $j$ of $n_l$ number of lines and $n_b$ number of buses, $\mathbf{\underline{Y}}^{\mathrm{Line}}= \begin{bmatrix} \mathbf{\underline{Y}}^{\mathrm{fr}}\\ 
\mathbf{\underline{Y}}^{\mathrm{to}} \end{bmatrix}$.}
\nomenclature[D, 23]{$\mathbf{E}^{max}$}{Vector of maximum energy capacity of $n_y$ storage devices and EVs (MWh).}%
\nomenclature[D, 24]{$\mathbf{\underline{S}}^{\mathrm{Line}}_\mathrm{max}$}{$\in \mathbb{C}^{2n_l \times 1}$ Vector of maximum rated capacity (MVA) of $2n_l$ Lines (from and to).}%

\nomenclature[E, 01]{$\mathbf{X}, \mathbf{x}_t$}{Set of all variables on the optimisation horizon, set of variables at one time step $t$.}

\nomenclature[E, 02]{$\boldsymbol{\mathcal{V}}_{t}, \boldsymbol{\Theta}_{t}$}{Voltage phase angles and magnitudes of $n_b$ number of buses at time $t$.}
\nomenclature[E, 03]{$v_{i,t}, \theta_{i,t}$}{Voltage phase angles and magnitudes of bus $i$, at time $t$.}

\nomenclature[E, 04]{$\boldsymbol{\mathcal{P}}_t^{\mathrm{g}}, \boldsymbol{\mathcal{Q}}_t^{\mathrm{g}}$}{$\in \mathbb{R}^{n_{g}\times 1}$ vector of active and reactive power of $n_g$ generators at time $t$.} 
\nomenclature[E, 05]{$p_{i,t}^{\mathrm{g}}, q_{i,t}^{\mathrm{g}}$}{$\in \mathbb{R}$ active and reactive power of generator located on bus $i$, at time $t$.}

\nomenclature[E, 06]{$\boldsymbol{\mathcal{SOC}}_{t}$}{State-of-charge of  $n_y$ storage devices and EVs at time $t$.}
\nomenclature[E, 07]{$soc_{i,t}$}{State-of-charge of storage device or EV located on bus $i$, at time $t$.}

\nomenclature[E, 08]{$\boldsymbol{\mathcal{P}}_{t}^\mathrm{ch}, \boldsymbol{\mathcal{P}}_{t}^\mathrm{dch}$}{Rate of charge and discharge active power of $n_y$ storage devices and EVs at time $t$.}
\nomenclature[E, 09]{$p_{i,t}^\mathrm{ch}, p_{i,t}^\mathrm{dch}$}{Rate of charge and discharge active power of storage device or EV located on bus $i$, at time $t$.}

\nomenclature[E, 10]{$\boldsymbol{\mathcal{Q}}_{t}^\mathrm{s}$}{Rate of reactive power provision of $n_y$ inverters of storage devices and EVs at time $t$.}
\nomenclature[E, 11]{$q_{i,t}^\mathrm{s}$}{Rate of reactive power provision of inverter of storage device or EV located on bus $i$, at time $t$.}
\nomenclature[F, 02]{$i, j, k$}{$i, j, k$ are applied in BATTPOWER context either as the index of bus number or generator number or storage device/EV.}
\nomenclature[F, 04]{$\mathrm{fr}, \mathrm{to}$}{from bus $i$ to $j$, to bus $i$ from $j$.}
\nomenclature[F, 06]{$ - $, $\sim$}{Overhead variable signs for linear and non-linear equations, for example $\widetilde{\mathbf{G}}(\mathbf{X})$ stands for vector of nonlinear equality constraints.}

\nomenclature[G, 01]{$i$}{Index of storage device/EV, transformer.}
\nomenclature[G, 02]{$j$}{Index of consumer fed by transformer $i$.}
\nomenclature[G, 03]{$t$}{Index of time.}
\nomenclature[G, 04]{$\kappa_i$}{Share of active load of $i^{th}$ transformer on total active loads of 32 MV-LV transformers.}
\nomenclature[G, 05]{$\phi_i$}{Share of reactive to active load of $i^{th}$ transformer.}
\nomenclature[G, 06]{$\Phi$}{Share of total reactive load on total active load of 32 MV-LV transformers.}
 \nomenclature[G, 07]{$\psi_{i,j}$}{Share of yearly energy consumption  of consumer $j$ fed by transformer $i$ on total consumers fed by transformer $i$.}
\nomenclature[G, 08]{$P^{PCC}_t, P^{gen}_t$}{Hourly active power generated of the main system feeder: PCC, and the second system feeder: gen.}
\nomenclature[G, 09]{$P^{tot}_{t}, Q^{tot}_{t}$}{Total active and reactive hourly power production of the entire system under study.}
\nomenclature[G, 10]{$P_{i,t}^{Trans}, Q_{i,t}^{Trans}$}{Hourly share of active and reactive power production on transformer $i$.}
\nomenclature[G, 11]{$P^d_{i,j,t}, Q^d_{i,j,t}$}{Hourly share of active and reactive power production on transformer $i$ and consumer $j$.}
\nomenclature[G, 12]{$f()$}{An index function which converts three-dimensional loads of consumer $j$ fed by transformer $i$ at the time $t$ to a two-dimensional array with $n_b$ number of buses and $T$ number of time steps refer to L.\ref{function1} and L.\ref{function2}.}
\printnomenclature

%% file: tex/appB.tex
\section{Mathematical Backbone of Power Flow} \label{appendixB}
The mathematical background to section \ref{problemFormulate} is elaborated here. Consider the vector of complex bus voltages in rectangular coordinates as illustrated  by  $\displaystyle\mathbf{\underline{V}}\in \mathbb{C}^{n_b\times 1}$, where $\mathbb{C}$ is a complex set. The voltage vector comprises complex elements as: $\mathrm{\underline{v}}_i  = \abs{\mathrm{v}_i}e^{j\mathrm{\theta}_i}$, where $\mathrm{\underline{v}}_i\in \mathbb{C}$, $\{v_i ,\theta_i \}\in  \mathbb{R}$ are the voltage magnitude and angle of the corresponding bus in polar coordinates, where $\mathbb{R}$ is a real set.  Moreover,  $\{\boldsymbol{\mathcal{V}}, \mathbf{\Theta}\} \in \mathbb{R}^{n_b \times 1}$ can be defined as vectors of real magnitude and angle of bus voltages.  In vector form, the relationship between rectangular and polar coordinates is shown as:
\begin{equation}
\label{eqn:Voltage}
\mathbf{\underline{V}}=\mathbf{diag}(\boldsymbol{\mathcal{V}})\ \exp(j\boldsymbol{\Theta})
\end{equation}
Line connectivity matrices of $\{\mathbf{C}^{\mathrm{fr}}, \mathbf{C}^{\mathrm{to}}\} \in \mathbb{B}^{n_l \times n_b}$ can be extracted from $\mathbf{BUS}^{\mathrm{from}}$ and $\mathbf{BUS}^{\mathrm{to}}$ vectors, such that $\mathrm{c}_{ik}^{\mathrm{fr}}=1$ if bus $k$ is connected to line $i$, and otherwise $\mathrm{c}^{\mathrm{fr}}_{ik}=0$, and the same holds for $\mathbf{C}^{\mathrm{to}}$. $\{\mathbf{\underline{V}}^{\mathrm{fr}}, \ \mathbf{\underline{V}}^{\mathrm{to}}\} \in \mathbb{C}^{n_l \times 1}$ are the vectors of complex bus voltages at line terminals, including ``from" and ``to" nodes, correspondingly. These vectors can be extracted using the connectivity matrices explained above shown in Eqs. \eqref{eqn:generalProblem1} and \eqref{eqn:generalProblem2}.
\begin{equation}
\label{eqn:generalProblem1}
\mathbf{\underline{V}}^{\mathrm{fr}} = \mathbf{C}^{\mathrm{fr}}\mathbf{\underline{V}}\\
\end{equation}
\begin{equation}
\label{eqn:generalProblem2}
\mathbf{\underline{V}}^{\mathrm{to}} = \mathbf{C}^{\mathrm{to}}\mathbf{\underline{V}}\\
\end{equation}
and therefore:
\begin{align}
\label{eqn:V_line}
\mathbf{\underline{V}}^{\mathrm{Line}}=
{\begin{bmatrix}
    \mathbf{\underline{V}}^{\mathrm{fr}} \\
    \mathbf{\underline{V}}^{\mathrm{to}} \\
\end{bmatrix}\ \mkern-10mu}_{2n_l \times 1} = 
{\overbrace{\begin{bmatrix}
    \mathbf{C}^{\mathrm{fr}} \\
    \mathbf{C}^{\mathrm{to}} \\
\end{bmatrix}}^{\mathbf{C}^{\mathrm{Line}}}\ \mkern-10mu}_{2n_l \times n_b} \mathbf{\underline{V}} 
\end{align}
In order to obtain the entire network flow, the vector of complex voltages $\mathbf{\underline{V}}$ has to be determined. This can be done using the well-known Kirchhoff's current law: the sum of external current injections at a bus $\mathbf{\underline{I}}^{\mathrm{bus}} \in \mathbb{C}^{n_b\times 1}$ is equal to the sum of internal - through lines - current injections to the same bus $\mathbf{\underline{I}}^{\mathrm{bus}}= \mathbf{\underline{Y}}^{\mathrm{bus}}\mathbf{\underline{V}}$, where
$\mathbf{\underline{Y}}^{\mathrm{bus}} \in \mathbb{C}^{n_b\times n_b}$ is the bus admittance matrix. The same principle is applied to compute the complex line current using complex bus voltages of line terminals, and line admittance matrix $\mathbf{\underline{Y}}^{\mathrm{Line}} \in  \mathbb{C}^{2n_l\times n_b}$. This is shown in (\ref{eqn:generalProblem3}) 
\begin{equation}
\label{eqn:generalProblem3}
\mathbf{\underline{I}}^{\mathrm{Line}}={\begin{bmatrix}
    \mathbf{\underline{I}}^{\mathrm{fr}} \\
    \mathbf{\underline{I}}^{\mathrm{to}} \\
\end{bmatrix}\ \mkern-10mu}_{2n_l \times 1} =
{\overbrace{\begin{bmatrix}
    \mathbf{\underline{Y}}^{\mathrm{fr}} \\
    \mathbf{\underline{Y}}^{\mathrm{to}} \\
\end{bmatrix}}^{\mathbf{\underline{Y}}^{\mathrm{Line}}}\ \mkern-10mu}_{2n_l \times n_b}
\mathbf{\underline{V}} 
\end{equation}
The relation between bus admittance and line admittance matrices is defined by (\ref{eqn:generalProblem4}).
\begin{equation}
\label{eqn:generalProblem4}
\mathbf{\underline{Y}}^{\mathrm{bus}} = (\mathbf{C}^{\mathrm{fr}})^\top \mathbf{\underline{Y}}^{\mathrm{fr}}+(\mathbf{C}^{\mathrm{to}})^\top \mathbf{\underline{Y}}^{\mathrm{to}}+\mathbf{\underline{Y}}^{\mathrm{shunt}}
\end{equation}
$\{\mathbf{\underline{Y}}^{\mathrm{fr}},  \mathbf{\underline{Y}}^{\mathrm{to}}\}\in \mathbb{C}^{n_l\times n_b}$, and
$\mathbf{\underline{Y}}^{\mathrm{shunt}} \in \mathbb{C}^{n_b\times n_b}$ is the matrix of shunt admittance. 
Finally, the external complex power injections into a bus $i$ can be computed as $\mathrm{\underline{s}}^{\mathrm{bus}}_i = \underline{\mathrm{v}}_i(\underline{\mathrm{i}}^{\mathrm{bus}}_i)^*$, whereas the  complex power flow over a line at the terminal $k$ can be calculated by  $\underline{\mathrm{s}}^{\mathrm{Line}}_k = (\mathbf{C}_k^{\mathrm{Line}}\mathbf{\underline{V}})(\underline{\mathrm{i}}^{\mathrm{Line}}_k)^*$, where $\{\underline{\mathrm{s}}^{\mathrm{bus}}_i ,  \underline{\mathrm{i}}^{\mathrm{bus}}_i , \underline{\mathrm{s}}^{\mathrm{Line}}_k, \underline{\mathrm{i}}^{\mathrm{Line}}_i  \} \in \mathbb{C}$ and $\mathbf{C}_k^{\mathrm{Line}} \in \mathbb{B}^{1 \times n_b}$  is the $k$\textsuperscript{th} element of $\mathbf{C}^{\mathrm{Line}}$ matrix. In summary, power injections into a bus and into a line can be extended in the form of vectors using (\ref{eqn:generalProblem5}).   
\begin{align}
\label{eqn:generalProblem5}
&\mathbf{\underline{S}}^{\mathrm{bus}} = \mathbf{diag}(\mathbf{\underline{V}})(\mathbf{\underline{I}}^{\mathrm{bus}})^* \ \ \in \ \ \mathbb{C}^{n_b \times 1}
\\
&\mathbf{\underline{S}}^{\mathrm{Line}} = \mathbf{diag}(\mathbf{\underline{V}}^{\mathrm{Line}})(\mathbf{\underline{I}}^{\mathrm{Line}})^* \ \ \in \ \ \mathbb{C}^{2n_l \times 1}
\end{align}

%% file: tex/appC.tex
\section{Size and Structure of Input Matrices}\ \label{appendixC}
BATTPOWER input matrices are introduced and elaborated in this section. Table \ref{tab:inputMatrix} summarises the input matrices fed into the BATTPOWER solver proposed in \cite{zaferanlouei_battpower_2021}
\begin{table*}[!htbp]
 \caption{Definition of Input Matrices }
\label{tab:inputMatrix}
\begin{threeparttable}
\begin{tabularx}{\linewidth}{ l c c X }
\toprule
&   \multicolumn{2}{c}{Size of Matrix}& \\
\cmidrule{2-3} 
Input     & $n$&  $m$ &Description \\ 
\midrule
$\mathbf{BUS}$     &$n_b$ &\tnote{1}& Examples can be found in \cite{zimmerman_matpower:_2011} \\ \midrule
$\mathbf{BRANCH}$  &$n_l$ &\tnote{1}&Examples can be found in \cite{zimmerman_matpower:_2011}\\ \midrule
$\mathbf{GEN}$     &$n_g$ &\tnote{1}&Examples can be found in \cite{zimmerman_matpower:_2011} \\ \midrule
$\mathbf{GENCOST}$ &$n_g$ &\tnote{1}&Examples can be found in \cite{zimmerman_matpower:_2011}\\ \midrule
$\mathbf{BATT}$    &$n_y$  &\tnote{1}& \texttt{BATT\_BUS}, \ \texttt{SOC\_OPT}, \ \texttt{PCH\_OPT}, \ \texttt{PDICH\_OPT}, \ \texttt{Q\_INJ\_OPT}, \  $\boldsymbol{\mathcal{SOC}}^\mathrm{max}$, \ $\boldsymbol{\mathcal{SOC}}^\mathrm{min}$, \ $(\boldsymbol{\mathcal{Q}}^\mathrm{s})^\mathrm{max}$, \ $(\boldsymbol{\mathcal{Q}}^\mathrm{s})^\mathrm{min}$, \ \texttt{MBASE}, \ $(\boldsymbol{\mathcal{P}}^\mathrm{ch})^\mathrm{max}$, \ $(\boldsymbol{\mathcal{P}}^\mathrm{dch})^\mathrm{max}$ \ \texttt{EFF\_CH} ($\mathbf{\Psi}^\mathrm{ch}$) \ \texttt{EFF\_DICH} ($\mathbf{\Psi}^\mathrm{dch}$) \\ 
\midrule
$\mathbf{AVBP}$   &$n_y$ &$T$&$\mathbf{AVBP}$ is a binary set ($\in \mathbb{B}^{n_y \times T}$) which is the availability matrix of active power provision of storage devices, such that  $\mathbf{AVBP}_{i,t}=1$ if the $i^{th}$ storage at $t^{th}$ time is available and connected to the grid, otherwise $\mathbf{AVBP}_{i,t}=0$, where $T$ is the optimisation horizon. \\ 
\midrule
$\mathbf{CONCH}$  &$n_y$ &$T$&$\mathbf{CONCH} \in \mathbb{B}^{n_y \times T}$ is the charge connectivity matrix in which $\mathbf{CONCH}_{i,t}=1$ if the $i^{th}$ storage at $t^{th}$ time has a charging option, otherwise $\mathbf{CONCH}_{i,t}=0$.  \\ \midrule
$\mathbf{CONDI}$  &$n_y$ &$T$&$\mathbf{CONDI} \in \mathbb{B}^{n_y \times T}$ is the discharge connectivity matrix such that $\mathbf{CONDI}_{i,t}=1$ if the $i^{th}$ storage at $t^{th}$ time has the available discharging option, otherwise $\mathbf{CONDI}_{i,t}=0$\tnote{2}.\\ \midrule
$\mathbf{AVBQ}$   &$n_y$ &$T$&$\mathbf{AVBQ} \in \mathbb{B}^{n_y \times T}$ is the availability matrix of reactive power provision of storage devices such that $\mathbf{AVBQ}_{i,t}=1$ if the $i^{th}$ storage at $t^{th}$ time has the available option for reactive power provision, otherwise $\mathbf{AVBQ}_{i,t}=0$. \\\midrule
$\mathbf{AVG}$    &$n_g$ &$T$&$\mathbf{AVG} \in \mathbb{B}^{n_g \times T}$ which is the availability matrix of generators within the optimisation time horizon and consequently $\mathbf{AVG}_{i,t}=1$ if the $i^{th}$ generator at $t^{th}$ time is available to inject power in the grid. \\\midrule
$\mathbf{SOCi}$   &$n_y$ &$T$&$\mathbf{SOCi} \in \mathbb{R}^{n_y \times T}$ is the matrix consisting of initial state of charge of $n_y$ storage devices over time $t \in \{1,...,T\}$.  A value for initial state of charge $\{0\leq \mathbf{SOCi}_{i,t}\leq 1\}$ is allocated for the $i^{th}$ storage device at time $t$ if and only if one of these conditions is satisfied: 1) $\mathbf{AVBP}_{i,t=1}=1$. 2) $\mathbf{AVBP}_{i,t-1}=0$ and $\mathbf{AVBP}_{i,t}=1$ (arrival definition), otherwise $\mathbf{SOCi}_{i,t}=0$.  \\ \midrule
$\mathbf{SOCMi}$  &$n_y$ &$T$&$\mathbf{SOCMi} \in \mathbb{R}^{n_y \times T}$ matrix which includes the minimum state of charge of $n_y$ storage devices through time $t \in \{1,...,T\}$. The state of charge of the $i^{th}$ storage device at the departure time of $t$ can be settled if one of these two conditions is satisfied: 1) $\mathbf{AVBP}_{i,t}=1$, $\mathbf{AVBP}_{i,t+1}=0$. 2) $\mathbf{AVBP}_{i,t=T}=1$. \\ \midrule
$\mathbf{PD}$ & $n_b$& $T$& Time series of active loads. \\ \midrule
$\mathbf{QD}$ & $n_b$& $T$& Time series of reactive loads.  \\ 
\bottomrule
\end{tabularx}
\begin{tablenotes}
\item[1] {User Defined}
\item[2] {Note that  $\mathbf{AVBP}_{i,t}=0$  means that the $i^{th}$ storage\textbackslash EV at time $t$ is not available; therefore, the same element in charge and discharge connectivity matrices must be zero: $\mathbf{CONCH}_{i,t}=0$ and $\mathbf{CONDI}_{i,t}=0$. The converse logic is not valid.}
\end{tablenotes}
\end{threeparttable}
\end{table*}

%% file: tex/appD.tex
\section{Input Data} \label{appendixD} 
In this appendix, we elaborate:
\begin{enumerate}[i.]
\item \textbf{Base-Load:} The  process of data preparation of hourly 856 consumers base-load.
\item \textbf{EV data:} The  process of data preparation of EVs charge profile, arrival and departure.
\end{enumerate}
\subsection{Estimation of Consumer's Base-Load} \label{estimate}
\begin{figure*}[!htbp]
\centering
\includegraphics[width=5.5 in , height=1 in]{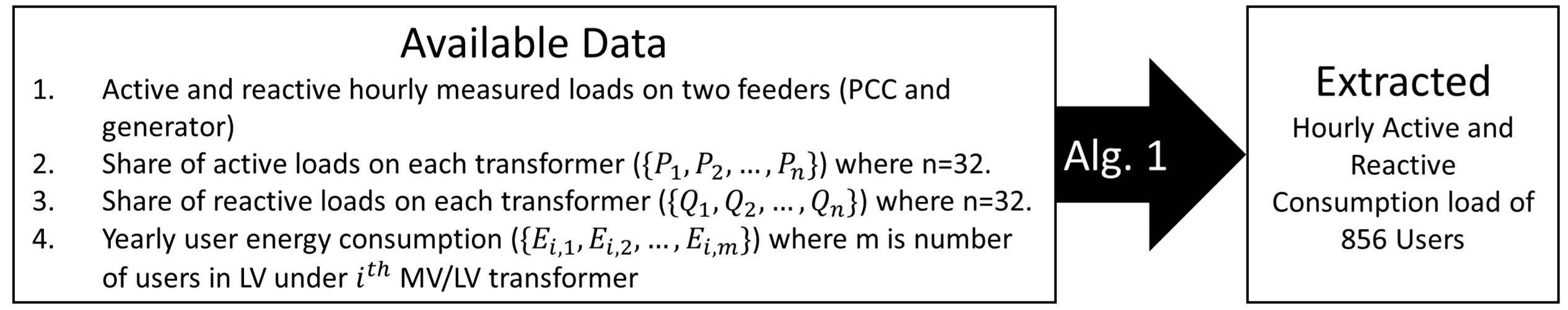}
\caption{Input and output of the Alg. \ref{Alg1} } 
\label{AlgorithmDes}
\end{figure*}

The hourly feeders (PCC and generator) load time series for 8208 hours (342 days) are imported from the local DSO, which is highly correlated to the ambient temperature \cite{lillebo_impact_2019, lillebo_impact_2018}. The feeder load data are from 25 Jan. 2012 until 31 Dec. 2012, which is assumed to be as the base-load with zero EV penetration profile (refer to subsection \ref{assumption}, item.\ref{zeroEV}). \\
From the input data available, shown in Fig. \ref{AlgorithmDes}, hourly active and reactive consumption loads of 856 consumers are estimated. 
First, four ratio factors are computed: 1) $\kappa_i$: share of active load of $i^{th}$ transformer on total active loads of 32 MV-LV transformers, L.\ref{kappa}. 2) $\phi_i$: share of reactive to the active load of $i^{th}$ transformer, 3) $\Phi$: share of total reactive load on total active load of 32 MV-LV transformers, and 4) $\psi_{i,j}$: share of yearly energy consumption  of consumer $j$ fed by transformer $i$ on total consumers fed by transformer $i$. Share of active loads on each transformer $\kappa_i$ is assumed to be constant throughout the simulation. The algorithm to compute consumers' hourly consumption for 8208 hours is depicted in Alg. \ref{Alg1}. 
\begin{algorithm}
\caption{Estimation of consumers' base-load}\label{Alg1}
\SetAlgoLined
\DontPrintSemicolon
$\forall i \in \{1,\dots,n\}$ \;
$\forall j \in \{1,\dots,m\}$\;
$\forall t \in \{1,\dots,T\}$\;
$\kappa_i = \frac{P_i}{P_{1}+P_{2}+\dots+P_{n}}$ \label{kappa}\;
$\phi_i = \frac{Q_{i}}{P_{i}} \quad \forall i \in {1,\dots,n}$\; \label{betai}
$\Phi = \frac{Q_{1}+Q_{2}+\dots+Q_{n}}{P_{1}+P_{2}+\dots+P_{n}}$\; \label{betan}
$\psi_{i,j} = \frac{E_{i,j}}{E_{i,1}+E_{i,2}+\dots+E_{i,m}}$ \;
\For{all Hours h}{
$P^{tot}_{t}=P^{PCC}_t+P^{gen}_t$\;
$Q^{tot}_{t}=P^{tot}_{t} \times \Phi$
}
\For{i= 1:all n Transformers}{
\For{t= 1:all HOURS h}{
$P_{i,t}^{Trans}=\kappa_i.P^{tot}_{t}$\;
$Q_{i,t}^{Trans}=\phi_i.P_{i,t}$
}
}
\For{i= 1:all n Transformers}{
\For{j= 1:all m consumers fed by $i^{th}$ transformer}{
\For{t= 1:all HOURS h}{
$P^d_{i,j,t}=\psi_{i,j}.P_{i,t}^{Trans}$\;
$Q^d_{i,j,t}=\psi_{i,j}.Q_{i,t}^{Trans}$
}
}
}
$\mathbf{PD}_{n_b \times T}=f({[P^d_{i,j,t}]}_{n \times m\times T})$\label{function1}\; 
$\mathbf{QD}_{n_b \times T}=f({[Q^d_{i,j,t}]}_{n \times m\times T})$ \label{function2}
\end{algorithm}
The EV optimisation horizon is selected to be 24 hours (96 time steps for 15-min resolution) from 12:00 PM until 12:00 PM the next day.
The highest base-load peak (system consumption) occurs 8:00 AM Feb. 3, 2012. The highest spot market price (NOK/MWh) belongs to 8:00 AM Feb. 2, 2012. 
\subsection{Charging profile data generation} \label{charge}
For the analysis of different scenarios, the presences of EVs with energy demand and corresponding charging power profiles relevant to their connection time needs to be provided. Therefore, the EV data must be generated with the available statistical data about driving patterns. The EV data generation method is detailed in the following subsections.
\subsubsection{EV drive distance}\label{distance}
Recent studies state the average distance driven in Norway is about 52 km \cite{bretteville-jensen_norwegian_2016, report_learning_2016}. The standard deviation of the drive distance specific to the selected locality is 22 km. The driving pattern and associated energy consumption depend on the traffic, weather conditions and associated road conditions. Daily EV energy demand for all EV owners (for the whole population) is generated by creating and normal distribution with a mean of 52 km and standard deviation of 22 km. A set of daily demand for a year is created for every EV with a fixed daily drive distance as the mean and 10\% of the fixed distance and standard deviation to accommodate the energy demand variation due to traffic, weather and road conditions. A summary of the data for EV charge profile generation is provided in Table \ref{tab:EV2}.
\subsubsection{EV distribution and energy consumption}\label{consump}
Reference \cite{elbilstatistikk} provides the statistics for EVs sold in Norway. The top 10 brands of EVs sold in Norway constitute 85\% of the total EV population. To simplify the modelling, the EV population is categorised based on their driving efficiency. 80\% of the EVs among the top brands, consumes in the range of 12 kWh to 18 kWh per 100 km. The rest of them consume 19 kWh to 25 kWh per 100 km. Therefore, the mean energy consumption for the two groups are assumed to be 17 kWh and 21 kWh. The same proportion of EVs is assumed in the simulated population.
\subsubsection{Arrival and departure time}\label{arridep}
The arrival and departure time of EVs predominantly depends on the weekly work time schedule. 78\% of the labour force in Norway has a regular shift work schedule which typically starts around 8:00 AM and ends around 18:00 PM \cite{sterud_working_nodate}.Total working hours per week are 40 hours. Though the daily working hours are 8:00 AM, they may begin between 8:00 AM and 10:00 AM, and end between 16:00 PM and 18:00 PM depending on the organisation type. The skilled labours and shift workers start their work schedules as early as 7:00 AM and end their work time as late 20:00 PM depending on the time they start. The average commutation time is one and half hours. The arrival time is generated with 17:00 PM as mean and 90 minutes as standard deviation for the given population. Though every EV owner has a fixed departure and arrival schedule, it is not precisely the same for every day. Therefore, 15 minutes' standard deviation is provided while the arrival time is generated for a year for every individual EV owner. The departure time is calculated by subtracting 9.5 hours from the arrival time.
\subsubsection{EV charging profile} \label{profile}
There are 2 types of domestic EV chargers commonly used in Norway. The power ratings are 2.3 kW (10 A) and 3.7 kW (16 A). There are very few EV owners, who have EVs with larger battery capacity use 11 kW (16 A 3 phase) charger. The percentages of presence of different chargers are given in Table \ref{tab:EV2}. The EVs which have the mean consumption of 17 kWh/100 km are charged with 2.3 kW chargers, and the rest are charged with 3.7 kW or 11 kW chargers.
The charging profiles are created with the assumptions that the EVs will start charging at the rated power of the chargers as soon as they arrive and the charging continues until the total drive demand for the day is fulfilled.
\begin{table}[htbp!]
\begin{center}
\scriptsize
\caption{Data for EV charge profile generation}
\label{tab:EV2}
\begin{threeparttable}
\begin{tabular}{l l} 
\toprule
 Mean daily drive distance &52 km  \\  
\midrule
Standard deviation of daily drive distance& 22 km\\
\midrule
Standard deviation of daily drive distance distribution&10\% \tnote{1}  \\
\midrule
Percentage of EV population that consume $\leq$ 18 kWh/100km &80\%\\
\midrule
Percentage of EV population that consume $\geq$ 18 kWh/100km &20\%\\
\midrule
Mean arrival time for the EV population& 17:00 hours\\
\midrule
Standard deviation of arrival time for the EV population& 90 min\\
\midrule
Standard deviation of daily arrival time for individual EV & 15 min\\
\midrule
Percentage of 230V, 10A chargers & 70\%\\
\midrule
Percentage of 230V, 16A chargers & 20\%\\
\hline
Percentage of 230V, 48A chargers & 10\%\\
\bottomrule
\end{tabular}
\begin{tablenotes}
\item[1] {of daily drive distance}
\end{tablenotes}
\end{threeparttable}
\end{center}
\end{table}

%% file: tex/appE.tex
\section{Market Data} \label{appendixE} 
\begin{figure*}[!htbp]
\centering
\includegraphics[width=3.5 in , height=2.2 in]{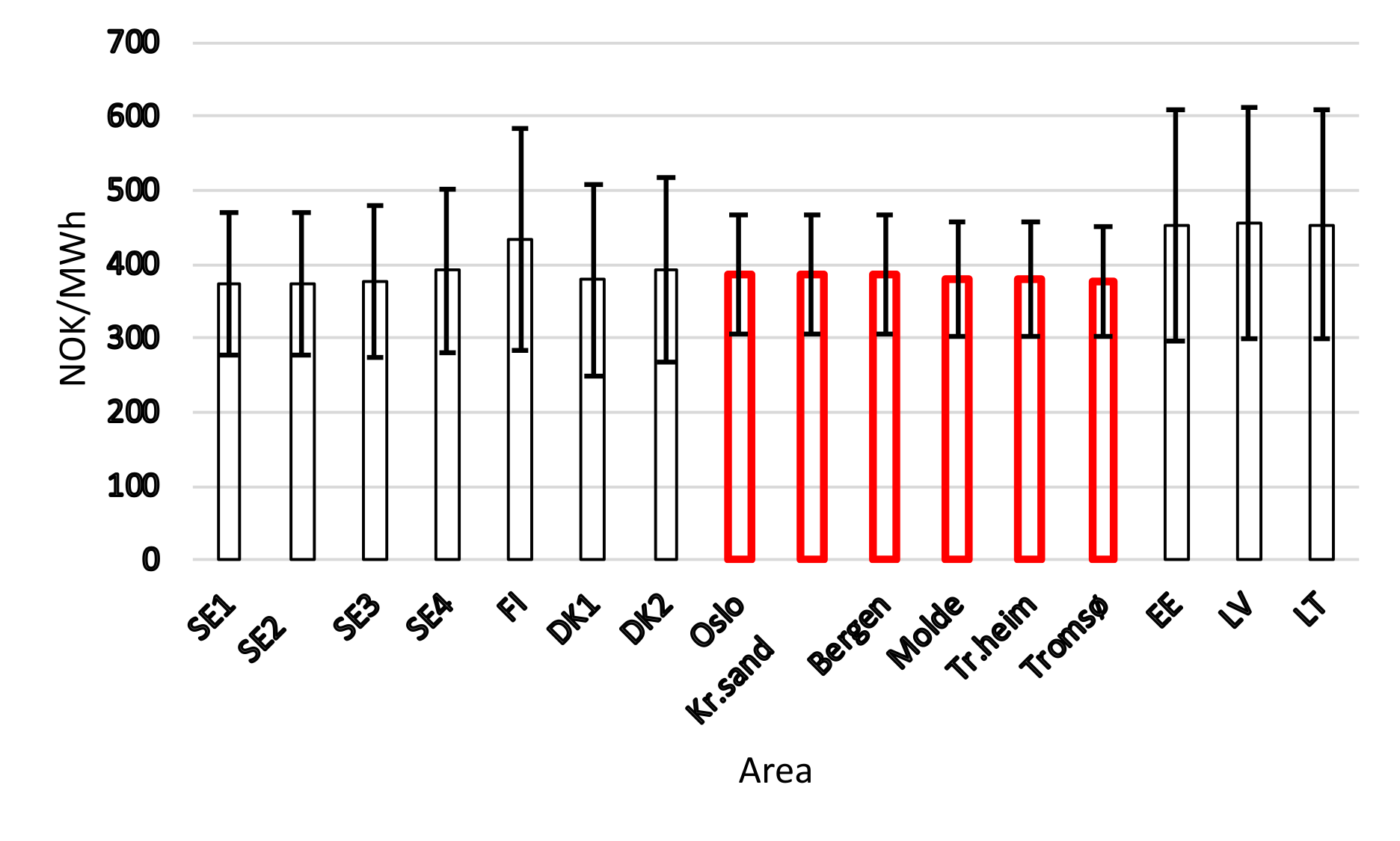}
\caption[]{Average and STD of 8760 hourly day-ahead prices (NOK/MWh) for the year of 2019. The average price of Norwegian areas are highlighted with red. They have lowest STD among other price areas.}
\label{fig:MarketData}
\end{figure*}

The goal of this appendix is to present that the Nowegian price areas mostly owns flat price profile in campare with other Nordic areas.\\
Table \ref{tab:EV} and Fig. \ref{fig:MarketData} show the average and standard deviation of spot price of the Nordic areas (NOK/MWh), year 2019. The Norwegian areas (Oslo, Kr.sand, Bergen, Molde, Tr.heim Tromsø) had the lowest STD. over the last year of 2019.\\
The average of price of Norwegian areas are highlighted in red to make them distiguishable from other area prices. The Norwegian price areas possess lower STD in compare with the other areas in Nordpool.
\begin{table*}[!htbp]
\footnotesize
\caption{Average and Standard deviation of 8760 hourly day-ahead prices (NOK/MWh) for the year 2019. Data can be found in \cite{nordpool}.}
\label{tab:EV}
\begin{tabularx}{\textwidth}{m s s s s s s s s s s s s s s s s} 
\toprule
 &{\scriptsize SE1}& {\scriptsize	SE2}	&{\scriptsize SE3}	&{\scriptsize SE4}&	{\scriptsize FI}& {\scriptsize	DK1}& {\scriptsize	DK2}&{\scriptsize	Oslo}&{\scriptsize	Kr.sand}	&{\scriptsize Bergen}	&{\scriptsize Molde}	&{\scriptsize Tr.heim}&	{\scriptsize Tromsø}&{\scriptsize	EE}	&{\scriptsize LV}&{\scriptsize 	LT}
 \\  
\midrule \midrule
AVE&373.7	&373.7	&377.9&	392.1&	434&	379.1&	392.4&	386.8&	386.6&	386.7&	379.6	&379.6	&377.3	&451.8&	455.9	&454.4
\\
STD&96.8&	96.8	&101.9	&111.3	&151.2&	129.6&	124.7	&81.2&	80.1	&80.5&	76.8	&76.8	&73.8&	156.1	&155.8&	155.7\\
\bottomrule 
\end{tabularx} 
\end{table*}

%% file: tex/appF.tex
\section{Danish Price Area as the Input Price Signal} \label{appendixF}
Analysis \ref{resultsanalysisc} is repeated with the same input data with the only expection of adopting DK2 price signal shown in Fig. \ref{Operational}, part (c). The purpose is to show with the high price volatility the savings increase using coordinated charge scheduling strategies. Note that i. the price DK2 is taken from one of the days with large price variations and ii. daily energy consumption and active loss are similar to the results presented in Table \ref{tab:termination}.
\begin{table*}[htbp!]
\begin{center}
\caption{[a) total energy production, b) active system loss, and c) system cost] in three different operational mode. Note that the Danish price profile of DK2 is the input data for optimisation here, and thus the obtained values in this table.}
\label{tab:danishanlysis}
\begin{threeparttable}
\begin{tabular}{l c c c c c} 
\toprule
 {\scriptsize Method} & \begin{tabular}[c]{@{}c@{}}{\scriptsize Daily\tnote{4} Energy}\\ {\scriptsize Consumption (MWh)}\end{tabular}  &  \begin{tabular}[c]{@{}c@{}}{\scriptsize Active Loss}\\ {\scriptsize (MWh)}\end{tabular}&  \begin{tabular}[c]{@{}c@{}}{\scriptsize System Cost}\\ {\scriptsize (NOK)}\end{tabular} &               \begin{tabular}[c]{@{}c@{}}{\scriptsize Daily Saving}\\ {\scriptsize (NOK)--(\%)}\end{tabular}  & \begin{tabular}[c]{@{}c@{}}{\scriptsize Yearly Saving}\\ {\scriptsize (NOK)}\end{tabular}  \\  
\midrule \midrule
 \tnote{1}&118.8 & 2.24 & 36,087.05 & -- & --\\
\midrule
 \tnote{2} & 118.767&2.18 & 33,143.26 & 2944 -- 8.2 \%  & 1,074,483\\
\midrule
 \tnote{3}&118.756&2.17 & 33,155.35 & 2932 -- 8.1\% & 1,070,071\\
\bottomrule
\end{tabular}
\begin{tablenotes}
{\footnotesize
\item[1]  Dumb Charging. 
\item[2]  MPOPF without operational limits. 
\item[3]  MPOPF with operational limits. 
\item[4]  Obtained based on the base-load input data of date of 12:00 PM 1 Feb. 2012 to 12:00 PM 2 Feb. 2012.}
\end{tablenotes}
\end{threeparttable}
\end{center}
\end{table*}

%% file: Main.bbl
\begin{thebibliography}{10}

\bibitem{EVoutlook}
IEA, ``Global ev outlook 2019.''
  \url{https://www.iea.org/reports/global-ev-outlook-2019}, 2019.

\bibitem{elbilstatistikk}
Elbilforeningen, ``Elbilstatistikk.'' \url{https://elbil.no/elbilstatistikk/},
  2020.

\bibitem{nveelbil}
NVE, P\"oyry, and DNV-GL, ``Kostnader i str\o mnettet - gevinster ved
  koordinert lading av elbiler.''
  \url{http://publikasjoner.nve.no/eksternrapport/2019/eksternrapport2019_51.pdf},
  2019.

\bibitem{carpentier1962contribution}
J.~Carpentier, ``Contribution a l’etude du dispatching economique,'' {\em
  Bulletin de la Societe Francaise des Electriciens}, vol.~3, no.~1,
  pp.~431--447, 1962.

\bibitem{chandy_simple_2010}
K.~M. Chandy, S.~H. Low, U.~Topcu, and H.~Xu, ``A simple optimal power flow
  model with energy storage,'' in {\em 49th {IEEE} {Conference} on {Decision}
  and {Control} ({CDC})}, pp.~1051--1057, Dec. 2010.
\newblock ISSN: 0743-1546.

\bibitem{wang_computational_2007}
H.~Wang, C.~E. Murillo-Sanchez, R.~D. Zimmerman, and R.~J. Thomas, ``On
  {Computational} {Issues} of {Market}-{Based} {Optimal} {Power} {Flow},'' {\em
  IEEE Transactions on Power Systems}, vol.~22, pp.~1185--1193, Aug. 2007.

\bibitem{wachter_implementation_2006}
A.~Wächter and L.~T. Biegler, ``On the implementation of an interior-point
  filter line-search algorithm for large-scale nonlinear programming,'' {\em
  Mathematical Programming}, vol.~106, pp.~25--57, Mar. 2006.

\bibitem{byrd_knitro_2006}
R.~H. Byrd, J.~Nocedal, and R.~A. Waltz, ``Knitro: {An} {Integrated} {Package}
  for {Nonlinear} {Optimization},'' in {\em Large-{Scale} {Nonlinear}
  {Optimization}} (G.~Di~Pillo and M.~Roma, eds.), Nonconvex {Optimization} and
  {Its} {Applications}, pp.~35--59, Boston, MA: Springer US, 2006.

\bibitem{kourounis_towards_2018}
D.~Kourounis, A.~Fuchs, and O.~Schenk, ``Towards the {Next} {Generation} of
  {Multiperiod} {Optimal} {Power} {Flow} {Solvers},'' {\em IEEE Transactions on
  Power Systems}, vol.~PP, no.~99, pp.~1--1, 2018.

\bibitem{sperstad_optimal_2016}
I.~B. Sperstad and H.~Marthinsen, ``Optimal power flow methods and their
  application to distribution systems with energy storage: a survey of
  available tools and methods,'' {\em report no. TR A7604, SINTEF Energy
  Research, Trondheim.}, 2016.

\bibitem{sperstad_energy_2019}
I.~B. Sperstad and M.~Korpås, ``Energy {Storage} {Scheduling} in
  {Distribution} {Systems} {Considering} {Wind} and {Photovoltaic} {Generation}
  {Uncertainties},'' {\em Energies}, vol.~12, p.~1231, Jan. 2019.

\bibitem{moghadasi_optimal_2016}
S.~Moghadasi and S.~Kamalasadan, ``Optimal {Fast} {Control} and {Scheduling} of
  {Power} {Distribution} {System} {Using} {Integrated} {Receding} {Horizon}
  {Control} and {Convex} {Conic} {Programming},'' {\em IEEE Transactions on
  Industry Applications}, vol.~52, pp.~2596--2606, May 2016.

\bibitem{capitanescu_experiments_2013}
F.~Capitanescu and L.~Wehenkel, ``Experiments with the interior-point method
  for solving large scale {Optimal} {Power} {Flow} problems,'' {\em Electric
  Power Systems Research}, vol.~95, pp.~276--283, Feb. 2013.

\bibitem{castillo2013computational}
A.~Castillo and R.~P. O’Neill, ``Computational performance of solution
  techniques applied to the acopf,'' {\em Federal Energy Regulatory Commission,
  Optimal Power Flow Paper}, vol.~5, 2013.

\bibitem{zaferanlouei_battpower_2021}
S.~Zaferanlouei, H.~Farahmand, V.~V. Vadlamudi, and M.~Korpas, ``{BATTPOWER}
  toolbox: Memory-efficient and high-performance {MultiPeriod} {AC} optimal
  power flow solver,'' pp.~1--1, 2021.
\newblock Conference Name: {IEEE} Transactions on Power Systems.

\bibitem{sojoudi_optimal_2011}
S.~Sojoudi and S.~H. Low, ``Optimal charging of plug-in hybrid electric
  vehicles in smart grids,'' in {\em 2011 {IEEE} {Power} and {Energy} {Society}
  {General} {Meeting}}, pp.~1--6, July 2011.
\newblock ISSN: 1944-9925.

\bibitem{masoum_smart_2011}
A.~Masoum, S.~Deilami, P.~Moses, M.~Masoum, and A.~Abu-Siada, ``Smart load
  management of plug-in electric vehicles in distribution and residential
  networks with charging stations for peak shaving and loss minimisation
  considering voltage regulation,'' {\em Transmission Distribution IET
  Generation}, vol.~5, pp.~877--888, Aug. 2011.

\bibitem{richardson_local_2012}
P.~Richardson, D.~Flynn, and A.~Keane, ``Local {Versus} {Centralized}
  {Charging} {Strategies} for {Electric} {Vehicles} in {Low} {Voltage}
  {Distribution} {Systems},'' {\em IEEE Transactions on Smart Grid}, vol.~3,
  pp.~1020--1028, June 2012.

\bibitem{chen_optimal_2012}
N.~Chen, T.~Q. Quek, and C.~W. Tan, ``Optimal charging of electric vehicles in
  smart grid: {Characterization} and valley-filling algorithms,'' in {\em 2012
  {IEEE} {Third} {International} {Conference} on {Smart} {Grid}
  {Communications} ({SmartGridComm})}, pp.~13--18, Nov. 2012.
\newblock ISSN: null.

\bibitem{oconnell_rolling_2014}
A.~O'Connell, D.~Flynn, and A.~Keane, ``Rolling {Multi}-{Period} {Optimization}
  to {Control} {Electric} {Vehicle} {Charging} in {Distribution} {Networks},''
  {\em IEEE Transactions on Power Systems}, vol.~29, pp.~340--348, Jan. 2014.

\bibitem{franco_mixed-integer_2015}
J.~F. Franco, M.~J. Rider, and R.~Romero, ``A {Mixed}-{Integer} {Linear}
  {Programming} {Model} for the {Electric} {Vehicle} {Charging} {Coordination}
  {Problem} in {Unbalanced} {Electrical} {Distribution} {Systems},'' {\em IEEE
  Transactions on Smart Grid}, vol.~6, pp.~2200--2210, Sept. 2015.

\bibitem{benetti_real-time_2015}
G.~Benetti, M.~Delfanti, T.~Facchinetti, D.~Falabretti, and M.~Merlo,
  ``Real-{Time} {Modeling} and {Control} of {Electric} {Vehicles} {Charging}
  {Processes},'' {\em IEEE Transactions on Smart Grid}, vol.~6, pp.~1375--1385,
  May 2015.

\bibitem{de_hoog_optimal_2015}
J.~de~Hoog, T.~Alpcan, M.~Brazil, D.~A. Thomas, and I.~Mareels, ``Optimal
  {Charging} of {Electric} {Vehicles} {Taking} {Distribution} {Network}
  {Constraints} {Into} {Account},'' {\em IEEE Transactions on Power Systems},
  vol.~30, pp.~365--375, Jan. 2015.

\bibitem{shao_layered_2015}
C.~Shao, X.~Wang, X.~Wang, and C.~Du, ``Layered and {Distributed} {Charge}
  {Load} {Dispatch} of {Considerable} {Electric} {Vehicles},'' {\em IEEE
  Transactions on Power Systems}, vol.~30, pp.~1858--1867, July 2015.

\bibitem{wang_integrated_2016}
D.~Wang, X.~Guan, J.~Wu, P.~Li, P.~Zan, and H.~Xu, ``Integrated {Energy}
  {Exchange} {Scheduling} for {Multimicrogrid} {System} {With} {Electric}
  {Vehicles},'' {\em IEEE Transactions on Smart Grid}, vol.~7, pp.~1762--1774,
  July 2016.

\bibitem{quiros-tortos_control_2016}
J.~Quirós-Tortós, L.~F. Ochoa, S.~W. Alnaser, and T.~Butler, ``Control of
  {EV} {Charging} {Points} for {Thermal} and {Voltage} {Management} of {LV}
  {Networks},'' {\em IEEE Transactions on Power Systems}, vol.~31,
  pp.~3028--3039, July 2016.

\bibitem{mehta_smart_2018}
R.~Mehta, D.~Srinivasan, A.~M. Khambadkone, J.~Yang, and A.~Trivedi, ``Smart
  {Charging} {Strategies} for {Optimal} {Integration} of {Plug}-{In} {Electric}
  {Vehicles} {Within} {Existing} {Distribution} {System} {Infrastructure},''
  {\em IEEE Transactions on Smart Grid}, vol.~9, pp.~299--312, Jan. 2018.

\bibitem{zhang_fast_2019}
J.~Zhang, M.~Cui, B.~Li, H.~Fang, and Y.~He, ``Fast {Solving} {Method} {Based}
  on {Linearized} {Equations} of {Branch} {Power} {Flow} for {Coordinated}
  {Charging} of {EVs} ({EVCC}),'' {\em IEEE Transactions on Vehicular
  Technology}, vol.~68, pp.~4404--4418, May 2019.

\bibitem{shi_model_2019}
Y.~Shi, H.~D. Tuan, A.~V. Savkin, T.~Q. Duong, and H.~V. Poor, ``Model
  {Predictive} {Control} for {Smart} {Grids} {With} {Multiple}
  {Electric}-{Vehicle} {Charging} {Stations},'' {\em IEEE Transactions on Smart
  Grid}, vol.~10, pp.~2127--2136, Mar. 2019.

\bibitem{kotsalos_horizon_2019}
K.~Kotsalos, I.~Miranda, N.~Silva, and H.~Leite, ``A {Horizon} {Optimization}
  {Control} {Framework} for the {Coordinated} {Operation} of {Multiple}
  {Distributed} {Energy} {Resources} in {Low} {Voltage} {Distribution}
  {Networks},'' {\em Energies}, vol.~12, p.~1182, Jan. 2019.

\bibitem{zaferanlouei_integration_2017}
S.~Zaferanlouei, M.~Korpås, H.~Farahmand, and V.~V. Vadlamudi, ``Integration
  of {PEV} and {PV} in {Norway} using multi-period {ACOPF} — {Case} study,''
  in {\em 2017 {IEEE} {Manchester} {PowerTech}}, pp.~1--6, June 2017.

\bibitem{cuffe_visualizing_2017}
P.~Cuffe and A.~Keane, ``Visualizing the {Electrical} {Structure} of {Power}
  {Systems},'' {\em IEEE Systems Journal}, vol.~11, pp.~1810--1821, Sept. 2017.

\bibitem{passengercars2012}
``Number of registered passenger cars in {Norway} from 2009 to 2020,
  https://www.statista.com/statistics/452433/norway-number-of-registered-passenger-cars/,''
  June 2021.

\bibitem{EVs2012}
``Elbilbestand, {Number} of electric cars and rechargeable hybrids in {Norway},
  https://elbil.no/elbilstatistikk/elbilbestand/,'' June 2021.

\bibitem{noauthor_statens_nodate}
``Statens vegvesen, statistics about vehicles in norway.''
  \url{https://www.vegvesen.no/}.

\bibitem{noauthor_tabell_nodate}
``Tabell 2 {Privathusholdninger} og personer per privathusholdning, etter
  fylke. 1960, 1970, 1980, 1990, 2001, 2011 og 2012.''
  \url{https://www.ssb.no/a/kortnavn/familie/tab-2013-01-17-02.html}.

\bibitem{nordpool}
``{Nord} {Pool}.'' \url{https://www.nordpoolgroup.com/}, 2020.

\bibitem{eriksen_rme_nodate}
A.~B. Eriksen and V.~Mook, ``{RME} {Rapport} 2/2020 {Proposed} changes to the
  design of network tariffs for low,'' {\em The Norwegian Energy Regulatory
  Authority (RME)}, pp.~1--21.

\bibitem{sortomme_optimal_2011}
E.~Sortomme and M.~A. El-Sharkawi, ``Optimal {Charging} {Strategies} for
  {Unidirectional} {Vehicle}-to-{Grid},'' {\em IEEE Transactions on Smart
  Grid}, vol.~2, pp.~131--138, Mar. 2011.

\bibitem{yao_hierarchical_2013}
W.~Yao, J.~Zhao, F.~Wen, Y.~Xue, and G.~Ledwich, ``A {Hierarchical}
  {Decomposition} {Approach} for {Coordinated} {Dispatch} of {Plug}-in
  {Electric} {Vehicles},'' {\em IEEE Transactions on Power Systems}, vol.~28,
  pp.~2768--2778, Aug. 2013.

\bibitem{ortega-vazquez_electric_2013}
M.~A. Ortega-Vazquez, F.~Bouffard, and V.~Silva, ``Electric {Vehicle}
  {Aggregator}/{System} {Operator} {Coordination} for {Charging} {Scheduling}
  and {Services} {Procurement},'' {\em IEEE Transactions on Power Systems},
  vol.~28, pp.~1806--1815, May 2013.

\bibitem{chung_electric_2019}
H.-M. Chung, W.-T. Li, C.~Yuen, C.-K. Wen, and N.~Crespi, ``Electric {Vehicle}
  {Charge} {Scheduling} {Mechanism} to {Maximize} {Cost} {Efficiency} and
  {User} {Convenience},'' {\em IEEE Transactions on Smart Grid}, vol.~10,
  pp.~3020--3030, May 2019.

\bibitem{zimmerman_matpower:_2011}
R.~D. Zimmerman, C.~E. Murillo-Sanchez, and R.~J. Thomas, ``{MATPOWER}:
  {Steady}-{State} {Operations}, {Planning}, and {Analysis} {Tools} for {Power}
  {Systems} {Research} and {Education},'' {\em IEEE Transactions on Power
  Systems}, vol.~26, pp.~12--19, Feb. 2011.

\bibitem{lillebo_impact_2019}
M.~Lillebo, S.~Zaferanlouei, A.~Zecchino, and H.~Farahmand, ``Impact of
  large-scale {EV} integration and fast chargers in a {Norwegian} {LV} grid,''
  {\em The Journal of Engineering}, vol.~2019, no.~18, pp.~5104--5108, 2019.

\bibitem{lillebo_impact_2018}
M.~Lillebo, ``Impact of {EV} {Integration} and {Fast} {Chargers} in a
  {Norwegian} {LV} {Grid} - {An} analysis based on data from a residential grid
  in {Steinkjer},'' 2018.

\bibitem{bretteville-jensen_norwegian_2016}
T.~Bretteville-Jensen, ``The {Norwegian} electric car controversy: {The}
  arguments and some empirical illustrations,'' {\em 60}, 2016.

\bibitem{report_learning_2016}
T.~Report, E.~Figenbaum, and M.~Kolbenstvedt, ``Learning from {Norwegian}
  {Battery} {Electric} and {Plug}-in {Hybrid} {Vehicle} users – {Results}
  from a survey of vehicle owners,'' p.~8, 2016.

\bibitem{sterud_working_nodate}
T.~Sterud, ``Working time in the {European} {Union}: {Norway}.''
  \url{https://www.eurofound.europa.eu/publications/report/2009/working-time-in-the-european-union-norway},
  2020.

\end{thebibliography}
